\documentclass[reqno,12pt]{amsart}
\usepackage{amssymb, amsmath, geometry}
\usepackage{graphicx}
\usepackage{amsmath,amssymb,amsthm,bm}
\usepackage{tikz}
\usepackage{geometry}
\usepackage{tikz-cd}

\usepackage{array}
\usepackage[colorlinks,linkcolor=blue!55!black,citecolor=blue!55!black]{hyperref}
\usetikzlibrary{arrows.meta,decorations.pathreplacing}
\usepackage{hyperref}
\hypersetup{
    colorlinks=true,
    linkcolor=blue,
}
\usepackage{blindtext}

\usepackage{etoolbox}

\makeatletter
\renewcommand*\l@section{%
  \@tocline{1}{0pt}{1em}{3em}{}}
\renewcommand*\l@subsection{%
  \@tocline{2}{0pt}{2em}{3em}{}}
\patchcmd{\@tocline}
  {\hfil}
  {\nobreak\leaders\hbox{\normalfont.\kern.4em}\hfill}
  {}
  {\PackageError{DMT}{Could not patch contents leaders}{}}
\makeatother

\usepackage{graphics}

\newcommand{\C}{\mathbb{C}}
\newcommand{\CP}{\mathbb{CP}}

\newcommand{\id}{\mathbb{I}}
\newcommand{\Id}{\mathbb{I}}
\newcommand{\zt}{\widetilde{\zeta}}
\newcommand{\tr}{\mathrm{Tr}}
\newcommand{\Ft}{\widetilde{F}}
\newcommand{\RP}{\mathbb{RP}}

\newcommand{\R}{\mathbb{R}}

\newcommand{\adj}{\mathrm{adj}}
\newcommand{\Z}{\mathbb{Z}}

\newcommand{\koniec}{\begin{flushright}  $\Box $ \end{flushright}}
\def\be{\begin{equation}}
\def\ee{\end{equation}}

\def\t{\tilde}
\def\p{\partial}
\def\tw{\widetilde{W}}
\def\tz{\widetilde{Z}}

\usepackage{bm}
\usepackage{enumerate}
\usepackage{xcolor}

\newtheorem{theo}{Theorem}[section] 
\newtheorem{prop}[theo]{Proposition}  
\newtheorem{lemma}[theo]{Lemma}
\newtheorem{defi}[theo]{Definition}
\newtheorem{rmk}[theo]{Remark}
\newtheorem{example}[theo]{Example}

\newcounter{mnotecount}[section]

\renewcommand{\themnotecount}{\thesection.\arabic{mnotecount}}

\newcommand{\mnote}[1]
{\protect{\stepcounter{mnotecount}}$^{\mbox{\footnotesize
$
\bullet$\themnotecount}}$ \marginpar{
\raggedright\tiny\em
$\!\!\!\!\!\!\,\bullet$\themnotecount: #1} }

\numberwithin{equation}{section}

\begin{document}
\date{30 October 2026}
\title{The twistor construction of toric gravitational instantons} 
\author{Maciej Dunajski}
\address{Department of Applied Mathematics and Theoretical Physics\\ 
University of Cambridge\\ Wilberforce Road, Cambridge CB3 0WA\\ UK.}
\email{m.dunajski@damtp.cam.ac.uk}
\author{Lionel Mason}
\address{The Mathematical Institute\\
Oxford University\\
Woodstock Road, Oxford OX2 6GG\\ UK.
}
\email{Lionel.Mason@maths.ox.ac.uk}
\author{Paul Tod}
\address{The Mathematical Institute\\
Oxford University\\
Woodstock Road, Oxford OX2 6GG\\ UK.
}
\email{paul.tod@sjc.ox.ac.uk}
\begin{abstract}
We give an explicit construction of all ALF and ALE toric gravitational instantons in terms of certain rational matrix functions of one complex variable on the Riemann sphere; these provide patching matrices for their corresponding  twistor holomorphic vector bundles. The reconstruction of the instantons arises from the solution to a Riemann--Hilbert splitting problem
and we show how this can be solved via linear algebra, resembling but not equivalent to the Belinski--Zakharov ansatz. This construction gives a complete classification 
together with an explicit form for all toric gravitational instantons
up to the resolutions of conical and orbifold singularities. We illustrate the formalism by computer--assisted numerical construction of  two  new ALF (but not AF) gravitational instantons with five nodes.
\end{abstract}
\maketitle

\tableofcontents

\section{Introduction}
Gravitational instantons are  non--compact, complete Riemannian
four--manifolds $(M, g)$ with
\be
\label{eq:0}
  \mathrm{Ric}_g=0,\qquad \int_M
  |\mathrm{Riem}_g|^2\,\mathrm{vol}<\infty.
\ee
The second condition in (\ref{eq:0}) implies that
gravitational instantons have curvature which is asymptotically decaying.
In the case when $(M, g)$ is hyper--K\"ahler, this condition led to a classification of instantons by the volume growth $\rho^d$, with the asymptotic flat models  $\R^d\times T^{4-d}, d=1, \dots, 4.$
The case $d=4$ is asymptotically locally Euclidean (ALE), $d=3$ is asymptotically locally flat (ALF), and $d<3$ correspond to ALG,  ALG*, ALH, and ALH* which has growth  $\rho^{4/3}$ \cite{ALreference}.
 
Our interest will be in the non--hyper--K\"ahler cases
where all known examples are ALF and
toric: they admit an isometric action of $T^2$ with at least one fixed point. A partial classification
of such solutions is based on the number $n\in
\mathbb{N}$ of nodes: the isolated fixed points of
the torus action. The Kerr and Taub--bolt instantons
have $n=2$,  and the Chen--Teo family 
\cite{CT1, CT2} has $n=3$. All these instantons
are Hermitian \cite{AA, TodHermitian}, and it has been proven \cite{BG, BGL} that 
the global Hermiticity condition implies $n\leq 3$.
It has recently been shown using existence theorems
for harmonic maps \cite{LiSun} that asymptotically flat (AF) non--Hermitian instantons exist for any $n>3$. An explicit example with $n=4$ was 
given by  Teo \cite{Teo}. It is not known whether ALF (but non AF)
non--hyper--K\"ahler instantons exist with $n>3$. It is also not known whether any non--hyper--K\"ahler ALE instantons exist.

Ricci--flat metrics
with two commuting Killing vectors have been extensively studied in General Relativity, where $(M, g)$ is Lorentzian and solutions have been constructed
using the Belinski--Zakharov dressing method \cite{BZ} for the
two--dimensional integrable harmonic map equation resulting from the Ricci--flat condition. This has proved powerful in the construction of toric gravitational instantons \cite{CT1,CT2,Teo}.  In \cite{W2, WM, FW, F} these integrable 
equations were incorporated into the twistor
integrability framework of \cite{MW}. 
The aim of our paper is to show that in Riemannian signature the twistor techniques lead
to explicit twistor description of all gravitational instantons in terms of their twistor bundles, which  in turn can be expressed uniquely up to certain choices in terms of their patching matrices.  The instanton Ernst potentials and metrics can be reconstructed from these $P(w)$ via a Riemann--Hilbert splitting in purely algebraic terms. The corresponding solutions of the harmonic map equation  (both the Gram matrix $J$ obtained directly from the metric, and the Ernst matrix $E$, see (\ref{eq:17})) can be read--off directly from the splitting. Unlike the Belinski--Zakharov method and its variants, the twistor approach is guaranteed to lead to all possible solutions.

The central object in our study is a symmetric two--by-two matrix $P=P(w)$ that we refer to as a twistor patching matrix, as it describes a holomorphic vector bundle over a reduced twistor space.  It
depends on one complex variable, is meromorphic with $n$ simple poles on the real axis, satisfies the algebraic condition $\mbox{det}P=-1$, and has asymptotic behaviour for large $w$
\[
{P_{\infty}}^{ALF} \sim\begin{pmatrix}
 1    & 0   \\
 0&  -1
\end{pmatrix}, \quad \mbox{or}\quad 
{P_{\infty}}^{ALE} \sim \begin{pmatrix}
\frac{1}{2w}     & 0   \\
 0&  -2w
\end{pmatrix}.
\]
We show in Theorem \ref{theopoles} that for any gravitational instanton it can be uniquely constructed from the analytic continuation of the Ernst potential on the axis near to $z=+\infty$. Let ${\mathcal C}_n$ be the set of all such matrices with ALF asymptotics.

The following theorem explicitly constructs
all elements of  ${\mathcal C}_n$ recursively 
(the ALE patching matrices
depend on $2n-2$ moduli, and are constructed similarly in Theorem \ref{theo:ALEm}.)
\begin{theo}
\label{theo_intro}
Let $P_n\in {\mathcal C}_n$ be regular at $w=a$  and let
\[
B=B(w)=\mathbb{I}+\frac{\mu}{w-a}uu^t\epsilon
\]
where $\mu\in \R$ is a constant,  $\epsilon=-\epsilon^t$ with $\epsilon_{01}=1$, and $u$ is a real column vector such that $u^t \epsilon P_n(a) \epsilon u=0$.
Then
\[
P_{n+1}=BP_nB^t \in {\mathcal C}_{n+1}
\]
and all elements of ${\mathcal C}_{n+1}$ arise from this construction starting from ${P_{\infty}}^{ALF}$.
The moduli space of elements of ${\mathcal C}_{n}$ is $2n-1$ dimensional, up to translation of $w$.
\end{theo}
Given a patching matrix $P$, a reconstruction of a toric Ricci--flat metric proceeds in the following steps:
\begin{enumerate}
\item Find invertible matrices $F,\Ft$ depending on $(r,z,\zeta \in \CP^1)$ such that
$F$ is holomorphic in $\zeta$ for $|\zeta|\leq1$,
$\Ft$ is holomorphic in $\zt=-\zeta^{-1}$ for $|\zt|\leq1$,
and
\[
   P(w)=F\cdot\Ft^{\,-1}, \quad\mbox{if}\quad
 w=z-\frac{r}{2}\Big(\zeta-\frac1\zeta\Big).
\]
\item Use the gauge freedom
$(F,\Ft)\mapsto(F\cdot C,\ \Ft\cdot C)$ with $C=C(r, z)$ to set 
$\Ft({\zt}=0,r, z)=\id$. Then
\[
  E(r,z)=F(\zeta=0,r,z)
\]
is the Ernst matrix (see (\ref{eq:17})) for a gravitational instanton. 
\item Conversely, given a toric Ricci--flat
metric with $J=J(r, z)$, the corresponding patching matrix is an analytic continuation
of $E(r=0, z)$.
\end{enumerate}
This Riemann--Hilbert splitting is solved in two ways in this paper. In Theorem \ref{prop:directsplit} we start with a partial fraction decomposition of $P$ into a sum of simple poles, and give a direct formula for $E$ and its splitting from the resulting data at the residues. Although the formulae underlying this approach bear a close resemblance to those used in the Belinski--Zakharov method, they differ in key details and our approach shows that all instantons can be obtained in this way.  
Alternatively, in Theorem
\ref{splitting_prop} we show that the Riemann--Hilbert problem underlying the splitting can also be solved recursively building on Theorem \ref{theo_intro}:
a splitting of $P_n\in \mathcal{C}_{n}$
gives rise to an explicit splitting of $P_{n+1}$.

The  instanton metric is recovered from $E$ by
solving a system of first--order linear PDEs
(\ref{eq:18}).  However, we also show that it can be obtained more directly algebraically. Its Gram matrix can be read off directly from the splitting matrices
as in Proposition \ref{prop:Jdirect}. Furthermore, we can construct the transverse conformal factor algebraically via a formula that is analogous to that arising in \cite{BZ,BM,KKV} in Proposition \ref{prop:Omega}.

The paper is organised as follows: in the next section we shall review the toric condition in Riemannian
Ricci--flat geometry, introduce the Yang equation for the matrix $J$, and the associated rod structure together with ALE and ALF examples. In \S\ref{twistorsec} we introduce the twistor patching matrix and directly analyze
all cases up to three nodes. We shall prove that these are exhausted by the Gibbons--Hawking metrics
in the self--dual (SD) or anti--self--dual (ASD) cases or Kerr--NUT if $n=2$. For $n=3$ we shall prove
\begin{theo}
\label{theo2_intro}
Let $P=P(w)$ be a patching matrix with $n=3$. The corresponding toric Ricci--flat metric belongs
to one of the three classes
\begin{enumerate}
\item Three--centre Gibbons--Hawking  \cite{GH78}.
\item ALE Euclidean Pleba\'nski--Demia\'nski metrics \cite{PD}, or their real analytic continuations.
\item ALF Chen--Teo metrics \cite{CT1,CT2}, or their real analytic continuations.
\end{enumerate}
All Chen--Teo and Pleba\'nski--Demia\'nski metrics arise from this construction.
\end{theo}
None of the ALE Pleba\'nski--Demia\'nski metrics are regular, but the Chen--Teo class contains a 2--parameter family of AF instantons.
Here a real analytic continuation is a member of the family with complex parameters for which the patching matrix, and hence the metric, is real (see the proofs of Propositions \ref{prop_PD} and \ref{propCT}). Some of the Chen--Teo continuations are positive definite. None of them is regular: they are Hermitian, with a Hermitian structure which is smooth wherever the metric is \cite{Der}, and so they are excluded by \cite{AL} (ALE) and \cite{BG} (ALF).
In \S\ref{section4}
we present a recursive construction of all 
ALE and ALF patching matrices, and
prove Theorem \ref{theo_intro} as well as its ALE counterpart. In \S\ref{sec:RH} we give two constructions of the Riemann--Hilbert holomorphic splitting of patching matrices: a direct procedure illustrated by the Kerr instanton in \S\ref{sec:dirsplit}
and a recursive procedure in \S\ref{sec:recsplit}. In 
\S\ref{sec:direct} we use the holomorphic splitting of $P$ to give a direct construction of the instanton metric, its Gram matrix $J$ and
the transverse conformal factor. Finally in \S\ref{sectionn5} we give computer--assisted evidence for two candidates for new ALF instantons with $n=5$ nodes.

In Appendix A we summarise the Ward twistor correspondence
for anti--self--dual Yang--Mills equations
and their reduction to an integrable harmonic map equation. In Appendices B and C we give some details of patching matrices
of the Pleba\'nski--Demia\'nski and the Chen--Teo instantons respectively.
\medskip

\noindent
{\bf Acknowledgements}
This work was initiated when the authors attended the 
programme Twistor Theory held at the Newton Institute in Cambridge in September 2024. It was completed during
the Simons Programme on Twistor Theory and its Applications held at the Institute of Mathematics of the Polish Academy of Sciences in Warsaw in September 2026. We are grateful
to the participants of both programmes, in particular to Lars Andersson, Bernardo Araneda, and James Lucietti,  for discussions.  
This research was supported by the Simons
Foundation grant SFI-MPS-T-Institutes-00010825, and by the State Treasury funds under the
project Organization of the Simons Semesters at the Banach Center - New Energies in 2026-2028
(MNiSW/2025/DAP/491).  LM is supported by the Simons Foundation Collaboration on Celestial Holography, STFC grant number ST/T000864/1 and the Leverhulme Research Fellowship 2026-561/9. PT acknowledges support from St John's College Oxford.

The construction of the $n=5$ ALF instantons
in \S\ref{sectionn5} was carried out with the assistance of Claude Opus 5.5 (Anthropic, model claude-opus-5-5).
The Maple worksheets and Python data files are available on request.  Chatgpt pro has also been used to check and suggest improvements to other parts of the paper.

\section{Toric gravitational instantons}
\label{section21}
Let $(M, g)$ be a Ricci--flat Riemannian manifold which admits
an effective isometric action of $T^2=U(1)\times U(1)$ with at least one fixed point. The corresponding 
two commuting Killing vectors $K_1, K_2$ generate the Lie algebra $\mathfrak{t}=\R^2$. The Lorentzian counterpart would be, for example, stationary and axisymmetric space--times. We shall see in 
\S\ref{sec:regularity} that in the Riemannian signature the  toric assumption puts more constraints on regularity.
\subsection{Normal form}
The Gram matrix $J_{AB}=g(K_A,K_B)$ with $A, B=1, 2$,  is a symmetric
matrix--valued function on the space of orbits $\Sigma=M/T^2$.  Set $r^2\equiv\det J$. The 
Ricci--flat equations imply that  $r$ is harmonic
on $\Sigma$, with harmonic conjugate $z$ 
such that $dr=*_\Sigma dz$. The function $z$ is defined up to an overall constant translation,
and $(r=+\sqrt{\mbox{det} J}, z)$ extend to coordinates $(r , z , \phi^1, \phi^2)$ on $M$ such that
\begin{equation}\label{eq:11}
  g=\Omega^2(dr^2+dz^2)+J_{AB}\,d\phi^A d\phi^B ,\qquad
  K_A=\frac{\partial}{\partial\phi^A},
\end{equation}
where $\Omega$ and $J$ depend on $(r\geq 0, z)$. These coordinates are global
\cite{KL2, AL}.
The remaining Ricci--flat equations consist
of the Yang equation
\begin{equation}\label{eq:12}
  r^{-1}\partial_r\!\left(r\,J^{-1}\partial_r J\right)
  +\partial_z\!\left(J^{-1}\partial_z J\right)=0
\end{equation}
together with a pair of linear PDEs
\begin{equation}\label{eq:13}
\begin{aligned}
 \partial_r\ln\Omega&=-\frac{1}{2r}
   +\frac{r}{8}\,\tr\!\Big(\big(J^{-1}\partial_rJ\big)^2-\big(J^{-1}\partial_zJ\big)^2\Big),\\[2pt]
 \partial_z\ln\Omega&=\frac{r}{4}\,\tr\!\Big(J^{-1}\partial_zJ\cdot J^{-1}\partial_rJ\Big).
\end{aligned}
\end{equation}
The integrability condition of \eqref{eq:13} is \eqref{eq:12}. Therefore, the objective is to 
solve \eqref{eq:12} and then obtain $\Omega$ from a quadrature.
\subsubsection{Integrability}
The equation
\eqref{eq:12} is a completely integrable system arising as a
symmetry reduction of the anti--self--dual Yang--Mills (ASDYM) equations \cite{W2} (see Appendix A).
It admits a Lax pair which it inherits {from} ASDYM: Given a matrix
$J=J(r, z)$,  define $Z=z+ir\in \mathbb{H}=\{(r, z)\in\mathbb{R}^2,
r>0\}$, and
\be
\label{laxpair}
L_+=D_+ +\frac{i}{\zeta+i} J^{-1}\p_Z J, 
\quad L_-=D_- -\frac{i}{\zeta-i} J^{-1}\p_{\bar{Z}} J
\ee
where
 $\zeta\in \mathbb{CP}^{1}$ is an auxiliary  complex parameter and
\begin{equation}
\label{eq:Dpm}
\begin{aligned}
D_+=\p_Z+\frac{i\zeta(\zeta-i)}{2r(\zeta+i)}\p_\zeta,\quad
D_-=\p_{\bar{Z}}-\frac{i\zeta(\zeta+i)}{2r(\zeta-i)}\p_\zeta.
\end{aligned}
\end{equation}
The condition $[L_+, L_-]=0$
is equivalent to the Yang equation
(\ref{eq:12}).
\subsection{Rod structure}
The torus action degenerates at the $r=0$
boundary of the orbit space 
$\Sigma=\mathbb{H}=\{(r, z), r>0\}$. 
Since $\mbox{det} J=0$ at this boundary, we have
$\mathrm{rank}\, J(0, z)<2$. Generically
this rank is equal to 1, and it vanishes
at $n$ isolated nodes $a_1<a_2<\dots<a_n$ where both Killing vectors are zero.
\begin{center}
\begin{tikzpicture}[scale=1.0]
  \draw[-{Latex[length=2mm]},thick] (-3.6,0)--(3.9,0) node[right]{$z$};
  \foreach \x/\lab in {-2/a_1,-1/a_2,1.2/a_{n-1},2.4/a_n}
    {\fill (\x,0) circle (1.9pt); \node[below] at (\x,-0.08) {\scriptsize $\lab$};}
  \node[above] at (-2.9,0.05) {\scriptsize $R_0$};
  \node[above] at (-1.5,0.05) {\scriptsize $R_1$};
  \node[above] at (0.2,0.05) {\scriptsize $\cdots$};
  \node[above] at (3.1,0.05) {\scriptsize $R_n$};
\end{tikzpicture}
\end{center}
The nodes divide the $\partial \Sigma$ boundary into
$n+1$ {rods}
\[
 R_0=(-\infty,a_1),\quad R_i=(a_i,a_{i+1}),\quad R_n=(a_n,\infty), \quad i=1, \dots, n-1,
\]
which in the Lorentzian setup correspond to
axes or black hole horizons \cite{Har, Hollands}. In the Riemannian case they can only be axes. 

The {rod vector} of $R_i$ is $l_i\in\R^2$ with $J(0,z)\,l_i=0$ for
$z\in R_i$ (and $l_i$ is constant in $z$).
To each rod we can adapt a basis
$(K, l_i)$ where the Killing vector $K$ is
non--zero on $R_i$. In this basis (see Lemma \ref{lemmaA11} for a proof of a stronger statement)
\begin{equation}\label{eq:14}
  J=\begin{pmatrix} j_i(z) & 0\\ 0& 0\end{pmatrix}+O(r^2),
  \qquad
  \Omega^2=\frac{c_i^2}{j_i}+O(r^2),
  \qquad j_i>0 \ \ \text{on }(a_i,a_{i+1}).
\end{equation}
\subsection{Regularity analysis}
\label{sec:regularity}
The rod formalism introduced so far can lead to
two kinds of singularity: conical defects
on the rods that prevent them from being smooth two--spheres, and orbifold singularities at the nodes.
The analysis of these singularities depends on which ones are removed first. Our choice (first conical, then orbifold) agrees with the  choice made in \cite{CT2} but is opposite to the choice of \cite{KL2, LiSun}.
\subsubsection{Non--conicality}
 Let $u=\p_{\phi}\in\ker\big(J(0,z)\big)$
on a rod $R$.  The normal form (\ref{eq:14})
implies that, near $r=0$,
\be
\label{metricneartherod}
  g\sim \Omega^2dr^2+J(u,u)\,d\phi^2
   =\Omega^2\Big(dr^2+{\tfrac{J(u,u)}{\Omega^2}}d\phi^2\Big),
   \qquad J(u, u)=O(r^2).
\ee
The surface gravity
\begin{equation}\label{eq:15}
  \kappa(u)^2\equiv\lim_{r\to0}\frac{J(u,u)}{\Omega^2 r^2}
\end{equation}
is constant along $R$
and satisfies $\kappa^2=c^{-2}$, where $c$ is as in \eqref{eq:14}.
The rescaled rod vector $l\equiv u/\kappa(u)$ is  $2\pi$--periodic. 
This eliminates the conical singularity on $R$.
Assume this has been done for all rods. The normalised rod vectors $l_0, l_1, \dots, l_n$ and the rod lengths
$(a_2-a_1), \dots, (a_n-a_{n-1})$ are the rod data of a toric instanton. Not all rod data are admissible---we shall look at this next.
\subsubsection{Non--orbifold conditions}
 Two commuting Killing vectors $K_1, K_2$ generate the Lie algebra
$\mathfrak{t}=\R^2$ and give rise to a lattice
$\Lambda=\{K\in\mathfrak t:\ \exp(2\pi K)=\mathrm{Id}\}\cong\Z^2$, 
and $T^2=\mathfrak t/\Lambda$.
Let us assume that there exist at least two nodes 
$(a_1, a_2)$
and pick $l_0, l_1$ as a basis of
$\Lambda$. Then
\[
  \begin{pmatrix} l_1\\ l_2\end{pmatrix}
  =\begin{pmatrix} 0&1\\ k& p_1\end{pmatrix}
   \begin{pmatrix} l_0\\ l_1\end{pmatrix},\qquad k,p_1\in\Z ,
\]
where $k=\pm1$. Replacing $l_2$ by $-l_2$ if necessary we may take
$k=-1$, so that $l_2=-l_0+p_1l_1$. In general one obtains the lattice conditions
\begin{equation}\label{eq:16}
   l_{i+1}=-l_{i-1}+p_i\,l_i ,\qquad p_i\in\Z .
\end{equation}
\medskip
\begin{center}
\begin{tikzpicture}[scale=0.95]
  \begin{scope}
    \draw[gray!25,very thin] (-0.4,-0.4) grid (4.4,3.4);
    \foreach \x in {0,...,4} \foreach \y in {0,...,3} \fill[gray!85] (\x,\y) circle (1.7pt);
    \fill[blue!15] (0,0)--(2,1)--(3,2)--(1,1)--cycle;
    \draw[blue!70,thick] (0,0)--(2,1)--(3,2)--(1,1)--cycle;
    \draw[-{Latex[length=2.4mm]},very thick] (0,0)--(2,1) node[below right]{$l_{i-1}$};
    \draw[-{Latex[length=2.4mm]},very thick] (0,0)--(1,1) node[above left]{$l_i$};
    \node[blue!70] at (3.55,1.35) {\footnotesize area $1$};
    \node at (2,-1.0) {\small (a) $|\det(l_{i-1},l_i)|=1$: smooth, $\C^2$};
  \end{scope}
  \begin{scope}[xshift=7.6cm]
    \draw[gray!25,very thin] (-0.4,-0.4) grid (4.4,3.4);
    \foreach \x in {0,...,4} \foreach \y in {0,...,3} \fill[gray!85] (\x,\y) circle (1.7pt);
    \fill[red!12] (0,0)--(2,1)--(3,3)--(1,2)--cycle;
    \draw[red!70,thick] (0,0)--(2,1)--(3,3)--(1,2)--cycle;
    \draw[-{Latex[length=2.4mm]},very thick] (0,0)--(2,1) node[below right]{$l_{i-1}$};
    \draw[-{Latex[length=2.4mm]},very thick] (0,0)--(1,2) node[above left]{$l_i$};
    \fill[red] (1,1) circle (2.6pt); \fill[red] (2,2) circle (2.6pt);
    \node[red!80] at (3.65,1.45) {\footnotesize area $3$};
    \node at (2,-1.0) {\small (b) $|\det|=3$: orbifold $\C^2/\Z_3$};
  \end{scope}
\end{tikzpicture}
\end{center}
If $|\det(l_{i-1},l_i)|=k$ then the half--open cell
spanned by $l_{i-1},l_i$ contains $k-1$ non--zero lattice points (shown in red);
$k=1$ is smooth, $k>1$ gives the orbifold singularity $\C^2/\Z_k$.
If the singularities are not present, then
each finite rod is a {bolt} $2$--sphere if the zeros of $j(z)$ are simple, and
$j_i'(a_i)=-j_i'(a_{i+1}) $ at the nodes.
\subsection{Ernst potential and patching matrix}
Let $J$ satisfy \eqref{eq:12}. Pick a rod $R$ and a Killing vector $K$ with
$K\neq0$ on $R$.
Ricci--flatness implies that there exists a twist potential $\psi:\Sigma\to\R$ defined up to  addition of a constant such that
$*(K\wedge dK)=d\psi$.
The matrix
\begin{equation}\label{eq:17}
  E=\frac{1}{|K|^2}\begin{pmatrix} 1 & -\psi\\ -\psi & \psi^2-|K|^4\end{pmatrix}
\end{equation}
can be shown to satisfy \eqref{eq:12} with $\det E=-1$. We shall
call $E$ the Ernst matrix associated to $g, R$ and $K$.
With $\omega\equiv
\dfrac{J_{12}}{J_{11}} $, the pair $J=(|K|^2,\omega)$ corresponds
to $E=(|K|^2, \psi)$ through
\begin{equation}\label{eq:18}
  \frac{r}{|K|^4}\,*_2 d\psi=d\omega .
\end{equation}
The central object in the twistor approach of 
\S\ref{twistorsec} will be the  matrix
\be
\label{tpmatrix}
P(z)=E(0,z).
\ee
The analytic continuation of $P(z)$ is the patching matrix for a holomorphic vector bundle over the
reduced twistor space (Appendix A). In \S\ref{twistorsec} it
 will be constructed from purely algebraic
data without the need for solving any differential equations. The Ernst matrix $E(r, z)$ is then recovered from the Riemann--Hilbert splitting of $P$
(again, an algebraic procedure which we will be able to carry out explicitly). Solving \eqref{eq:18} gives $J$, and finally the quadrature \eqref{eq:13} gives 
$\Omega$.
\begin{example}[Euclidean Schwarzschild Instanton]
This example, and the resulting form of $P(z)$
illustrates all features of the {\em forward} construction of the patching matrix starting from $J$.
The metric is
\[
  g=\Big(1-\frac{2m}{\rho}\Big)d\tau^2+\Big(1-\frac{2m}{\rho}\Big)^{-1}d\rho^2
    +\rho^2\big(d\theta^2+\sin^2\theta\,d\phi^2\big),\qquad m>0,\ \rho\geq 2m,
\]
with $U(1)\times SU(2)$ isometry, so $T^2$ is generated by
$K_1=\partial/\partial\tau$, $K_2=\partial/\partial\phi$. Then, w.l.o.g.
\[
  J=\begin{pmatrix} 1-\dfrac{2m}{\rho} & 0\\[4pt] 0 & \rho^2\sin^2\theta\end{pmatrix},
  \qquad
  r=\sqrt{\rho(\rho-2m)}\,\sin\theta,\qquad z=(\rho-m)\cos\theta .
\]
There are three rods, $R_0:\theta=\pi$, $R_1:\rho=2m$ (the bolt $S^2$),
$R_2:\theta=0$, with nodes at $z=-m$ and $z=m$.
The normalised rod vectors at $R_0$ and $R_2$
are $l_0=\p_\phi$ and $l_2=-\p_\phi$. The non--normalised
rod vector at $R_1$ is $\p_\tau$. Using 
\eqref{eq:15} 
we find $\kappa(\partial_\tau)={(4m)^{-1}}$. Hence the rod
vectors are
\[
  l_0=\partial_\phi,\qquad l_1=4m\,\partial_\tau,\qquad
  l_2=-\partial_\phi=-l_0+0\cdot l_1
\]
which gives the lattice condition \eqref{eq:16} with
 $p_1=0$. There $\tau\sim\tau+8\pi m$, and
the period is the inverse Hawking temperature. The topology is $M=\R^2\times S^2$ which illustrates
the general feature: the Euler characteristic of $M$
is equal to the number of nodes.

To compute $P(z)$
take $R=R_2$, the
semi--infinite rod $z\in(m,\infty)$, and $K=\partial/\partial\tau$. Then
$K\wedge dK=0$, so $\psi=0$ and
\[
  E=\begin{pmatrix} |K|^{-2} & 0\\ 0 & -|K|^{2}\end{pmatrix},
  \qquad\text{so}\qquad
  P(z)\equiv E(r=0,z)=
  \begin{pmatrix} \dfrac{z+m}{z-m} & 0\\[6pt] 0 & -\dfrac{z-m}{z+m}\end{pmatrix}.
\]
The matrix $P$ is rational with simple poles and rank--one residues
\[
  2m\begin{pmatrix}1\\0\end{pmatrix}\otimes(1,0)\ \text{ at } z=m,
  \qquad
  2m\begin{pmatrix}0\\1\end{pmatrix}\otimes(0,1)\ \text{ at } z=-m .
\]
\end{example}
\subsection{ALE and ALF}
\begin{defi}
\label{defi:ALF}
A gravitational instanton is 
\begin{itemize}
\item {Asymptotically locally flat} (ALF) if it
approaches an $S^1$--bundle over $\R^3$, or the quotient of such a bundle by a finite group acting freely, at infinity.  
\item {Asymptotically locally Euclidean} (ALE) if it approaches $\R^4/\Gamma$
 at infinity, where $\Gamma$ is a discrete
subgroup of $SO(4)$.
\end{itemize}
\end{defi}
\noindent
In the ALF case, there exist asymptotic
coordinates $(\rho, \tau, \theta, \phi)$ with
\[
g\sim\big(d\tau+2N\cos\theta\,d\phi\big)^2+d\rho^2
+\rho^2\big(d\theta^2+\sin^2\theta\,d\phi^2\big).
\]
If $N=0$ the bundle is trivial, $g\sim S^1\times\R^3$, and ALF reduces to AF (asymptotically flat); the Schwarzschild and Kerr
instantons are examples.  In the toric case $K_1=\partial_\tau$, $K_2=\partial_\phi$ and
\[
  \text{AF}:\ l_0=\pm l_n, \qquad\qquad
  \text{ALF}:\ |\det(l_0,l_n)|=4|N|\neq0.
\]
ALF instantons have volume growth $\sim\rho^3$.
For a quotient end the coordinates $(\tau,\phi)$ are not separately periodic: they are identified by the lattice $\Lambda$ of \S\ref{sec:regularity}, and the circle generated by the primitive element of $\Lambda$ along $\p_\tau$ can have non--trivial isotropy on the semi--infinite rods (see \S\ref{sec:n5results}).
\vskip3pt
An example of an ALE metric is $\R^4$. It is AE with $\Gamma=\{\mathrm{Id}\}$,
\[
g=d\rho^2+\tfrac14\rho^2(\sigma_1^2+\sigma_2^2+\sigma_3^2)
\]
where $\sigma_1, \sigma_2, \sigma_3$ are left--invariant one--forms on $SU(2)$.
Kronheimer \cite{K1} showed that for any $\Gamma$ (cyclic $A_N$, dihedral $D_N$,
tetrahedral, octahedral, and icosahedral) there exists
an ALE anti--self--dual gravitational instanton. 
It is not known \cite{Nak} whether there exist ALE Ricci--flat metrics that are neither self--dual nor anti--self--dual.
ALE instantons have volume growth $\sim\rho^4$.
\subsection{Gibbons--Hawking metrics}
Explicit classes of hyper--K\"ahler ALE and ALF toric instantons
arise from an axisymmetric harmonic function
$V=V(r, z)$ and the Gibbons--Hawking ansatz:
\begin{equation}\label{eq:21}
  g=V\big(dr^2+r^2d\phi^2+dz^2\big)+V^{-1}\big(d\tau+A\big)^2,
  \qquad dV=*_{\R^3}dA.
\end{equation}
The  commuting Killing vectors are
$K_1=\partial/\partial\tau, K_2=\partial/\partial\phi$,
and regular metrics are given by the
multi--centre ansatz
\begin{equation}\label{eq:22}
   V=\varepsilon+\sum_{k=1}^{n}\frac{m_k}{\big(r^2+(z-a_k)^2\big)^{1/2}},
   \qquad \varepsilon\in\{0,1\},
\end{equation}
with the nodes $a_1<a_2<\dots<a_n$ on the $z$--axis. The regularity of the metric
requires all mass points to be equal: $m_1=m_2=\dots=m_n=m$. For large $\rho=\sqrt{r^2+z^2}$ one has
\begin{itemize}\itemsep1pt
 \item $\varepsilon=1$: \emph{ALF};\quad
       $g\sim\big(dr^2+r^2d\phi^2+dz^2\big)+\big(d\tau+A\big)^2$.
 \item $\varepsilon=0$: \emph{ALE};\quad
       $\displaystyle g\sim\frac{nm}{\rho}\Big(d\rho^2+\rho^2\big(d\theta^2+\sin^2\theta\,d\phi^2\big)\Big)
        +\frac{\rho}{nm}\big(d\tau+A\big)^2 .$
       Setting $\rho=\hat\rho^{\,2}/(4nm)$ and $\tau=2nm\,{\hat\psi}$ gives
       $g\to d\hat\rho^{\,2}+\tfrac14\hat\rho^{\,2}(\sigma_1^2+\sigma_2^2+\sigma_3^2)$
       on $\R^4/\Z_n$, the $A_{n-1}$ ALE space.
\end{itemize}
\paragraph{Gibbons--Hawking catalogue.}
\begin{center}
\begin{tabular}{c|l|l}
 & ALE & ALF\\
 $n=1$ & $\R^4$ (AE) & SD/ASD Taub--NUT\\
 $n=2$ & Eguchi--Hanson & multi--centre \\
 $n>2$ & multi--centre & multi--centre. 
\end{tabular}
\end{center}
The rod structure can be computed directly
in $(r, z)$ coordinates as the Gibbons--Hawking
metric is in the canonical form \eqref{eq:11}.
One has
$\ker\big(J(0,z)\big)\sim\partial_\phi-W\partial_\tau$ 
where $A=W d\phi$, and $W=W(z)$ can be found directly. This gives the lattice condition \eqref{eq:16}
\[
   l_{i+1}=-l_{i-1}+2\,l_i\qquad(\text{all }p_i=2)
   \]
or equivalently $l_{i+1}=-i\,l_0+(i+1)\,l_1$ with every cell in the fan of area $1$
\medskip
\begin{center}
\begin{tikzpicture}[scale=1.15]
  \draw[gray!25,very thin] (-0.4,-0.4) grid (2.4,4.4);
  \foreach \x in {0,1,2} \foreach \y in {0,...,4} \fill[gray!85] (\x,\y) circle (1.5pt);
  \foreach \y in {0,...,3}{
    \fill[blue!12] (0,0)--(1,\y)--(2,{\y+1})--(1,{\y+1})--cycle;
    \draw[blue!45] (0,0)--(1,\y)--(2,{\y+1})--(1,{\y+1})--cycle;}
  \foreach \y in {0,...,4} \draw[-{Latex[length=2mm]},very thick] (0,0)--(1,\y);
  \node[right] at (1,0) {$l_0$};
  \node[right] at (1,1) {$l_1$};
  \node[right] at (1,2) {$l_2$};
  \node[right] at (1,3) {$l_3$};
  \node[right] at (1,4) {$l_4$};
\end{tikzpicture}
\end{center}
\vspace{2pt}
\medskip
To compute the Ernst data \eqref{eq:17}
take $K=\partial/\partial\tau$ and find 
$\psi=V^{-1}$ and,
\[
E=\begin{pmatrix} V & -1\\ -1& 0\end{pmatrix},\qquad \det E=-1 .
\]
Therefore the patching matrix for the Gibbons--Hawking class adapted to the external rod $R_n$ is 
an analytic continuation of
\be
\label{GHpatch}
P_{GH}=\begin{pmatrix} V(z) & -1\\ -1& 0\end{pmatrix}, \quad V(z)=\varepsilon+\sum_{k=1}^n\frac{m_k}{z-a_k}.
\ee
\subsection{Asymptotics of $P$ on the external rod}
The ALF asymptotics can be modelled on $\R^3\times S^1$, or the Euclidean Schwarzschild which gives
\be
\label{ALFasymp}
  \begin{pmatrix} \dfrac{z+m}{z-m} & 0\\[6pt] 0 & -\dfrac{z-m}{z+m}\end{pmatrix}
  \ \xrightarrow[\ z\to\infty\ ]{}\
  g_0=\begin{pmatrix}1&0\\0&-1\end{pmatrix}
  \quad\text{---\ the asymptotic model for ALF.}
\ee
To find the ALE asymptotics take the flat metric on $\R^4$ given by $g=d\rho_1^2+d\rho_2^2+
J_{AB}d\phi^Ad\phi^B$ with $J=\mbox{diag}(\rho_1^2, \rho_2^2)$. Following the algorithm in
\S\ref{section21} we find $z=\frac{1}{2}(\rho_1^2-\rho_2^2), r=\rho_1\rho_2$ so  on the rod with $\rho_2=0$ 
\be
\label{ALEasymp}
  \begin{pmatrix}\frac{1}{2z}&0\\0&-2z\end{pmatrix}
  \quad\text{---\ the asymptotic model for ALE.}
\ee
Equivalent (up to a conjugation \eqref{C_transform}) asymptotic models arise from the Gibbons--Hawking patching matrix (\ref{GHpatch}) with $V=1$ for ALF and $V=z^{-1}$ for ALE.
\section{Twistor patching matrix}
\label{twistorsec}
The Lax formulation \eqref{laxpair} of the Yang equation \eqref{eq:12}
allows one to read off a solution once the kernel
of the Lax pair is known. The twistor approach gives
a direct route to this kernel. The Riemann--Hilbert
splitting of $P$ (see \S\ref{sec:RH}) gives rise to
an invertible matrix $F=F(z, r, \zeta)$ holomorphic for $|\zeta|<1$ and such that $L_{\pm} F^{-1}=0$.
In fact $E(r, z)$ can be read--off from $F$ without the need for integration,
and the same goes for $J$
(see \S\ref{sec:direct1}).
We shall summarise the procedure
below, and give details in subsequent sections with the general theory in Appendix A.
\begin{itemize}
\item
The forward construction of $P$ from the $J$--matrix 
is as follows
\[
  J\ {\xrightarrow{\ \mathrm{Ernst}\ }}\ E
    \ \longrightarrow\ P(z)=E(0,z)
    \ \xrightarrow{\ \text{analytic continuation}\ }\ P(w),
\]
where 
\begin{equation}\label{eq:31}
  w=z-\frac{r}{2}\Big(\zeta-\frac1\zeta\Big)
  \end{equation}
is in the kernel of operators $D_{\pm}$ given by
(\ref{eq:Dpm}). The matrix
 $P(w)$ has determinant $-1$, is rational in $w$ with simple real poles, and on a semi--infinite rod
 its asymptotics is \eqref{ALFasymp} and \eqref{ALEasymp}.
 \item
The converse of this will be  achieved in two steps
\begin{enumerate}\itemsep2pt
\item Find (algebraically, with no differential equations) the most general
$P(w)$ with
  \begin{itemize}\itemsep0pt
   \item $n$ simple poles,
   \item asymptotics \eqref{ALFasymp} or \eqref{ALEasymp}.
   \item $\det P=-1$.
  \end{itemize}
\item Find invertible matrices $F,\Ft$ depending on $(r,z,\zeta)$ such that
  \begin{itemize}\itemsep0pt
   \item $F$ is holomorphic in $\zeta$ for $|\zeta|\leq1$,
   \item $\Ft$ is holomorphic in $\zt=-\zeta^{-1}$ for $|\zt|\leq1$,
  \end{itemize}
  and
\begin{equation}\label{eq:33}
   P(\zeta,r,z)=F\cdot\Ft^{\,-1}.
\end{equation}
\end{enumerate}
This holomorphic splitting is a non--abelian Riemann--Hilbert problem \cite{W1, W2, WM}. The
solution is unique up to the gauge freedom
$(F,\Ft)\mapsto(F\cdot C,\ \Ft\cdot C)$ with $C=C(r, z)$. This gauge freedom can be used to set 
$\Ft({\zt}=0,r, z)=\id$. Then
\begin{equation}\label{eq:34}
  E(r,z)=F(\zeta=0,r,z)
\end{equation}
satisfies \eqref{eq:12} with $\det E=-1$, and all  solutions to the Yang equation
arise this way \cite{WM}.
\end{itemize}
Although the first part of this follows from the general theory of \cite{W2, WM} as summarized in Appendix A, the identification $P(z)=E(0,z)$ gives distinct  $P(z)$ for different choices of rod as $E(0,z)$ is not in general analytic through a node. We fix this ambiguity by taking $P(z)$ to be the analytic continuation of the $E(0,z)$ restricted to the positive semi-infinite rod.

\subsection{Analytical properties of the patching matrix}\label{meromorphy}
In this section we shall establish the analytical properties of the patching matrix which result from the regularity of $(M, g)$.
\begin{theo}
  \label{theopoles}
  The twistor patching matrix for toric ALF or ALE instantons
  is meromorphic with simple poles on the real axis.
 \end{theo}

\noindent
{\bf Proof.} This follows from the correspondence of \cite{W1,W2,WM} as described in Appendix A and below, followed by some detailed examination of  the behaviour at nodes following both from \cite{F,FW} and further analysis below.  

We first review the relevant material from \cite{W2,WM} that gives the analyticity of $P(w)$ on the complement inside $\CP^1$ of the real axis $\RP^1$ together with the chosen semi--infinite rod going out to $+\infty$.

\smallskip

\noindent
{\bf Review of reduced twistor correspondence:} This follows on  from Appendix A but should be self contained.  After symmetry reduction the twistor correspondence reduces to a geometric correspondence between complexified reduced space-time
$\mathbb{C}^2$ with coordinates $(r, z)$ and a reduced twistor space $\mathcal{R}$ via the double fibration \cite{WM} 
\[
\begin{tikzcd}[column sep=small]
 & \mathcal{F} \arrow[dl, "p_1"'] \arrow[dr, "p_2"] & \\
 \mathbb{C}^2 & & \mathcal{R}.
\end{tikzcd}
\]
The three-dimensional correspondence space $\mathcal{F} = \mathbb{C}^2 \times \mathbb{CP}^1$ has coordinates
$(r, z, \zeta)$, and the fibration over $\mathbb{C}^2$ is
$p_1(r, z, \zeta) = (r, z)$. The $\mathbb{CP}^1$ fibres of $p_1$
are Riemann spheres with coordinate $\zeta$, the $\zeta$-spheres. The other fibration is generated by the vector fields $D_\pm$ \eqref{eq:Dpm} and leads to the projection
$p_2(r, z, \zeta) = w$, given by
\begin{equation}\label{eq:incidence}
  w = z + \tfrac{1}{2}r\left(\zeta^{-1} - \zeta\right),
\end{equation}
where $w$ is a local coordinate on $\mathcal{R}$ which is a covering of the Riemann sphere with affine coordinate $w$. 
To be precise about this structure, we must fix the domain
 $U\subset \C^2$ on which we are to define the  Ward twistor transform. 
 Then, this covering is obtained as the space of leaves of the fibration $p_2$ restricted to $U\times \CP^1$. The fibres of $\mathcal{R}$ over a point $w\in \CP^1$ are given by the leaves
 in $\mathcal{F}$ given by 
 \begin{equation}
     \zeta=\frac{z-w \pm\sqrt{(z-w)^2 + r^2}}{r}
 \end{equation} 
holding  $w$ constant. The leaf is connected if the square root has a branchpoint in $U$, whereas it has two components  elsewhere.  When $U$ is connected and simply connected and contains a connected piece $I$ of the axis, the map from $\mathcal{R}$ to the 
$w$-sphere is $1:1$ on  
\begin{equation}
    \tilde U=\{w=z\pm ir,  (r,z)\in U\}\, ,
\end{equation}
which is also connected and simply connected, and otherwise it is $2:1$. Thus $\mathcal{R}$ is two copies of the $w$ Riemann sphere glued together on $\tilde U$ by the identity map on $w$.

With this geometry, the Lax pair  \eqref{laxpair} lives on the correspondence space $\mathcal{F}$ acting on a trivial rank two vector bundle ${\mathcal E}'\rightarrow \mathcal{F}$.  The Ward bundle 
${\mathcal E}\rightarrow \mathcal{R}$ is defined so that its fibre at a point of $\mathcal{R}$ is the vector space of solutions to the Lax pair over the corresponding fibre of $p_2$.  It is shown in \cite{WM} that the regularity of the Ernst matrix at the interval on the axis that we have assumed implies that the restriction of ${\mathcal E}$ to each of the Riemann spheres making up ${\mathcal E}$ is trivial.
This is a special case of the axis-simple description of \cite{WM} with $p=q=0$ so that the bundle is trivial over each copy of the $w$-sphere.

We can now define $P$ to be the holomorphic matrix function on $\tilde U$ that glues  the  trivializations on each $w$-sphere. 
This determines the solution on $U$ as follows. For fixed
$r$ and $z$, \eqref{eq:incidence} is a $2:1$ map from the
$\zeta$-sphere to the $w$-sphere, branched at $\zeta = \pm i$, that is, over $w = z \pm ir$. The Ward bundle $\mathcal{E}\rightarrow \mathcal{R}$, is that descended from twistor space and the analogue of Ward's `triviality on twistor lines'
condition from Theorem \ref{theoWard}  is that the pullback $p_2^*\mathcal{E}$ is
trivial on the fibres of $p_1$, i.e.\ on the $\zeta$-spheres. To reconstruct the Ernst matrix on $U$,  pull $P$ back to the
$\zeta$-spheres, and split as
\begin{equation}\label{eq:splitting}
  P\big(w(r, z, \zeta)\big) = F(r, z, \zeta)\,\widetilde F(r, z, \zeta)^{-1}.
\end{equation}
The $F,\tilde F$ provide the transformation from the global frames on the two $w$-spheres to that on the $\zeta$-sphere.  The Ernst matrix is then recovered as 
\begin{equation}
    E= F(r, z, 0)\tilde F^{-1}(r, z, \infty)\, .
\end{equation}
This is the splitting that underlies the paper.  It is easily seen from this that in the limit as $r\rightarrow 0$ $P, F ,\tilde F$ become independent of $\zeta$ and this equation identifies $P(z)=E(0,z)$.


\smallskip
{\bf Meromorphic extension:} This follows from the work of Fletcher and Woodhouse \cite{F,FW}, in particular \S4.4, Lemma 4.1 of \cite{F}, who studied the behaviour of $P(w)$ in a neighbourhood of the intersection of the horizon and axis.  Under the sign flip of $\det J$, this gives the Euclidean signature neighbourhood of a node. Focusing on a node at $z=0$, there is a patching matrix on each side, say $P^+$ and $P^-$ and they prove, by a local space-time analysis, that up to conjugation by constants, we must have    
\begin{equation}
P^+=\begin{pmatrix}
{\alpha}/{2w}&\beta\\ \beta& 2w\gamma
\end{pmatrix}, \quad z> 0\, , \qquad 
P^-=\begin{pmatrix}
2w\alpha&\beta\\\beta& {\gamma}/{2w}
\end{pmatrix}  \, , \quad z<0, \label{Fletcher}
\end{equation}
where $\alpha, \beta$ and $\gamma$ are analytic in $w$ near $0$, with $\alpha\gamma-\beta^2=-1$ so that $\det P^\pm=-1$. 
Up to these constant conjugations $P^-=D\,P^+D$ with $D={\rm diag}\big(2(w-a_i),\,(2(w-a_i))^{-1}\big)$, so the continuation of $P^+$ through $w=a_i$ determines $P^-$; this provides a meromorphic extension of $P^-$ through to the next interval and so-on by induction in such a way that, apart from possibly at $\infty$, no new poles appear.   

\smallskip

{\bf Remaining lemmas:} We wish to check that we do not find worse than the simple poles that we expect at the nodes and that the simple poles we do have do indeed result in a smooth four--manifold.  This will follow from substituting normal forms
of the $J$--matrix and the conformal factor near a rod into the metric,
and requiring that this metric restricts to a regular metric on the bolt
two--sphere. These normal forms will be found in Lemma \ref{lemmaA11}. We will
also need to show that the end--points cannot be removed from the rod
as it would lead  to manifolds with cylindrical ends.
\vskip 5pt
Start from the canonical form of the metric (\ref{eq:11})
where the Gram matrix $J$ and the conformal factor $\Omega$ depend on $(r, z)$  and satisfy 
(\ref{eq:12}) and (\ref{eq:13}).
The twistor patching matrix is then the analytic continuation of the Ernst potential (\ref{eq:17}).
For this to have a double pole, the squared norm $|K|^2$
would have to have a double zero. We shall show that this can not happen.
\begin{lemma}
\label{lemmaA11}
Let $J$ satisfy the Yang equation  (\ref{eq:12}) with $\mbox{det}{(J)}=r^2$.
There exists a basis of Killing vectors, such that in a neighbourhood of a
rod $R_i\equiv (a_{i}, a_{i+1})$
\be
\label{near_rod}
J(r, z)=\begin{pmatrix}
j(z) & 0   \\
0 & 0
\end{pmatrix}
+O(r^2), \quad \Omega^2(r, z)=\frac{c^2}{j(z)}+O(r^2)
\ee
where $j>0$ when $z\in (a_{i}, a_{i+1})$ and $c$ is a constant.
\end{lemma}
\noindent
{\bf Proof.}
Choose a basis of the torus $(\p/\p\phi^1, \p/\p\phi^2)$ such that
$\p/\p\phi^2$ vanishes on the rod $R_i$, and $j\equiv g(\p/\p\phi^1, \p/\p\phi^1)|_{R_i}$
is a positive function of $z$ in the interior of the rod.
Therefore\footnote{We could just quote (\ref{near_rod})  if it was the case that $J$ can only
smoothly depend on $(z, r^2)$. This is known to be true in the Lorentzian black hole context \cite{weinstein}, but not in
 Riemannian.
If one assumes that  Weinstein's result \cite{weinstein} holds in Riemannian,
then this Lemma is in \cite{KL}.}
\be
\label{intermediaJ}
  J(r, z)=\begin{pmatrix}
j+j_1 r+j_2 r^2 & k_1r+k_2 r^2   \\
 k_1r+k_2 r^2 &   j^{-1}(1+{k_1}^2) r^2+l_3 r^3 
\end{pmatrix}
+
\begin{pmatrix}
O(r^3) & O(r^3)   \\
 O(r^3) &   O(r^4) 
\end{pmatrix},
\ee
with $j>0, j_1, j_2, k_1, k_2, l_3$ functions of $z$,
is the  general ansatz
which guarantees $\det{(J)}=r^2+O(r^3)$. The residue of the pole at $r^{-1}$ 
in the $12$--component of the
$J$--equation vanishes iff $k_1=0$. Then the $11$--component
implies   $j_1=0$. Finally the $22$--component forces $l_3=0$.

Now  expand the equations for the conformal factor
(\ref{eq:13}) in $r$. This yields, using $l_3=0$,
\[
  \p_z\ln{\Omega}=-\frac{1}{2}\p_z\ln{j}+O(r^2), \quad \p_r\ln{\Omega}=O(r)
\]
which gives (\ref{near_rod}).
\koniec

{\bf Regularity of resulting four--manifold.}  We wish to show that the instanton $M$ is smooth at the nodes arising from simple poles (subject to conicality conditions on transverse conformal factor).  
The metric near the rod $R_i$ takes the form
\[
  g=[j^{-1}(c^2 dr^2+r^2(d\phi^2)^2)]+ [c^2 j^{-1}  dz^2+j(d\phi^1)^2]+\dots
\]
where $\dots$ include the terms of order $O(r^2)$ as well as the cross terms
$d\phi^1d\phi^2$. If we assume  that $\phi^2$ is
normalised with period $2\pi$, then $c^2=1$ to avoid the conical singularity at
the rod. We now want the second factor
\be
\label{hinrod}
h\equiv j^{-1}dz^2+j(d\phi^1)^2
\ee
to be a smooth metric on a two--sphere. This is the case
(see, e.g. Lemma 4.1 in \cite{CDKL}) iff $j$ has simple zeros
at the end points $(a_{i}, a_{i+1})$, and  $j'(a_{i})=-j'(a_{i+1})$
and (to avoid the conical singularities
on the sphere) $\phi^1$ is periodic with the period $4\pi/|j'(a_{i})|$.

If one of the end points (say $a_{i}$) is $z=0$, then near this end--point
$j=\mbox{const}\times z$, and setting $z=\hat{z}^2$ puts $h$ in a locally flat form,
which will be free of conical singularities if the periodicity of $\phi^1$ is
adjusted. To finish--off the proof we need to rule out non--spherical rods.
If $a_{i}=0$ is a double zero of $j(z)$, i. e. $j(z)=\gamma z^2+O(z^3)$
then, near the turning point, (\ref{hinrod}) is a constant curvature hyperbolic metric 
$\gamma^{-1} z^{-2} dz^2+\gamma z^2 d(\phi^1)^2$ with a hyperbolic  end.
But this implies that the Ricci flat metric $g$ has two ends
which, in the ALE and ALF cases, is forbidden by the 
Cheeger--Gromoll splitting theorem \cite{CG}: If $(M, g)$ is complete
then it contains a geodesic (which turns out to be a line) connecting any two points on opposite ends, but then
$g$ is a product
metric\footnote{This argument does not apply to Einstein--Maxwell instantons.
  In fact the extreme 
  Reissner--Nordstr\"om instanton does have two ends and the function $j(z)$
  has a double zero with the metric on $M$ approaching
  $\mathbb{H}^2\times S^2$.}
  on $\R\times \Sigma$ which it can not be, at least
with the cubic or quartic growth needed for ALF or ALE.
\koniec
This theorem allows us to make a general ansatz
for the $P$--matrix. In the ALF case it takes the form
\begin{equation}
\label{Pal}
P^{ALF}=g_0+\sum_{i=1}^n\frac{e_i\, {v_i {v_i}^t}}{w-a_i}\, , \qquad \det P=-1, \quad\mbox{where}\quad g_0=
\begin{pmatrix}
1 & 0   \\
0 & -1
\end{pmatrix}, \quad e_i=\pm 1
\end{equation}
for some constant $2$--vectors $v_i$
and in the ALE case
\begin{equation}
P^{ALE}=\begin{pmatrix}
0 & 0   \\
0 & -2
\end{pmatrix}
(w-a_0)+\sum_{i=1}^n\frac{e_i\, {v_i {v_i}^t}}{w-a_i}\, , \qquad \det P=-1. \label{P-ansatze-ALE}
\end{equation} 
The ALE asymptotic form \eqref{ALEasymp} is then recovered using the translational freedom in $w$ to set $a_0=0$.
What needs a justification is the rank--one condition
of the residues of $P$. This follows from imposing
the determinant condition $\mbox{det}P=-1$
and removing the singularities at double poles.
\subsection{Solving the determinant condition directly}
Given the statement of Theorem \ref{theopoles}
as well as the asymptotic properties of $P$ we can set
\be
\label{PnewQRS}
P=
\frac{1}{\Delta}\begin{pmatrix}
R & Q\\
Q & S
\end{pmatrix}
\ee
with
\[
\Delta=(w-a_1)(w-a_2)\cdots(w-a_n)=w^n+\alpha_{n-2}w^{n-2}+\dots+\alpha_0
\]
where the freedom in translating $w$ has been used to set $\sum_{i=1}^{n}a_i=0$. The coefficients of the polynomials $(R, Q, S)$ are restricted by the determinant condition
\be
\label{detP}
\det{P}=-1, \quad\mbox{or}\quad RS-Q^2+\Delta^2=0.
\ee
There is a freedom in $P$
given by \cite{F, FW}
\be
\label{C_transform}
P\rightarrow CPC^{{t}}
\ee
where $C$  is of the form
\be
\label{C_transform2}
C=
\begin{pmatrix}
\delta & 0\\
\alpha+\beta w & \delta^{-1}
\end{pmatrix}
\ee
where $\alpha$ has the effect of adding a constant to the twist potential and $\beta$ corresponds to the freedom in choosing a Killing vector which does not vanish on the given rod.
\subsubsection{ALE}
For ALE the large $w$ asymptotics is
 $P\sim\mbox{diag}(\frac{1}{2w}, -2w)$, so we set
\be
\label{ALEcoeffs}
R=\frac{1}{2}w^{n-1}+\sum_{i=0}^{n-2}r_i w^i, \quad
Q=2Nw^{n-1}-\frac{L}{2}w^{n-2}+
\sum_{i=0}^{n-3}q_iw^i, \quad
S=-2w^{n+1}+\sum_{i=0}^n s_i w^i.
\ee
Applying (\ref{C_transform}) with $\delta=1, \alpha=-4N, \beta=0$ can be used to set $N=0$.
This leaves $4n-2$ parameters $(a, r, q, s, L)$ and $2n$ equations resulting from (\ref{detP}), as
the coefficient of $w^{2n}$ in the determinant vanishes identically. We solve the top
$2n-1$ equations corresponding to the coefficients
of $w^{2n-1}, \dots, w$ for $2n-1$ parameters
$(s, q)$. This can be done without branching. The final condition resulting from evaluating the determinant at $w=0$ gives a polynomial equation of
degree $2n-2$ for $L$. This factorises into two equations of degree $n-1$. The moduli space has dimension
\be
\label{NALE}
\mathcal{N}_{ALE}=2n-2.
\ee
We will examine the first three cases:
$n=1$ corresponds to the flat Euclidean space $\R^4$. $n=2$ is necessarily ASD and comes from Gibbons--Hawking. The solution is double Eguchi--Hanson metric with different masses which add up to a fixed constant, and depends on the distance between two centres, and the difference between the masses. It is free of conical singularities iff the masses are equal, which gives the one--parameter family of Eguchi--Hanson instantons. 
To see it set $a_1=a$ so that $\Delta=w^2-a^2$. The four equations resulting from (\ref{detP}) can be solved
\[
s_0= 8r_0(2{r_0}^2-a^2), \quad 
s_1=-4(2{r_0}^2-a^2), \quad {s_2}=4r_0, \quad
L=\pm(2a^2-8{r_0}^2)
\]
substituting this with the $-$ sign into $P$, and performing the transformation (\ref{C_transform})
with
\[
C=\begin{pmatrix}
{2} & 0\\
 2 r_0-w & \frac{1}{2}
\end{pmatrix}
\]
yields the Gibbons--Hawking patching matrix (\ref{GHpatch})
with the potential given by 
\[
V=\frac{m_+}{w-a}+\frac{m_-}{w+a}, \quad
\mbox{with}\quad m_{\pm}=1\pm\frac{2r_0}{a}.
\]
This solution is free of conical singularities iff $r_0=0$
in which case it is the Eguchi--Hanson instanton.
\begin{prop}
\label{prop_PD}
The patching matrix for $n=3$ ALE is equivalent, up to (\ref{C_transform}), to 
the multi Eguchi--Hanson with three different masses, or the four--parameter family of Pleba{\'n}ski--Demia{\'n}ski metrics, or its real analytic continuation. 
\end{prop}
\noindent
{\bf Proof.}
Five equations in (\ref{detP})
can be solved for
$\{s_0, s_1, s_2, s_3, \alpha_0\}$. The final equation for $L$ is a quartic which factorises
into a quadratic and two linear forms. The roots are 
\[
L_{ASD}=\pm(-8 r_{1}^{2}-2 \alpha_{1}+4 r_{0}),\quad
L_{PD}=-\frac{2 \left(-4 r_{1}^{3}-\alpha_{1} r_{1}+6 r_{0} r_{1}+\alpha_{0}\right)}{\sqrt{r_{1}^{2}-2 r_{0}}}.
\]
Taking $L=L_{ASD}$ with the $-$ sign and considering (\ref{C_transform}) with
\[
C=\begin{pmatrix}
\frac{2}{\beta} & 0\\
2\beta r_1-\beta w & \frac{\beta}{2}
\end{pmatrix}
\]
gives the general $n=3$ ALE Gibbons--Hawking
potential (\ref{GHpatch}) where, by a choice of 
$\beta$
we can set $m_1+m_2+m_3=1$. 

To complete the proof of Proposition \ref{prop_PD} we need to show that taking the second root $L=L_{PD}$ yields a patching matrix equivalent to that which results from the general Pleba\'nski--Demia\'nski metric. This is done as follows. In  Appendix B  we consider the Pleba\'nski--Demia\'nski metric and construct its patching matrix. It depends on four parameters
$(x_1, x_2, x_3, \gamma)$ and a matrix 
$C$ can be found which brings the patching matrix to the canonical form (\ref{ALEcoeffs}). The top left component is sufficient to read--off the ALE moduli
$(r_0, r_1, \alpha_0, \alpha_1)$ as functions of $(x_1, x_2, x_3, \gamma)$. We then compute the Jacobian
\[
dr_0\wedge dr_1\wedge d\alpha_0\wedge d\alpha_1=
f(x_i, \gamma) dx_1\wedge dx_2\wedge dx_3\wedge d\gamma
\]
and find that $f\neq 0$ unless some of the $x_i$s are equal. Therefore the transformation can be inverted
to express $(x, \gamma)$ as functions of $(r_0, r_1, \alpha_0, \alpha_1)$,
and so each non--ASD $n=3$ ALE patching matrix
results from some Pleba\'nski--{Demia\'nski} solution, possibly with complex parameters.

The Jacobian shows that the map from the parameters to the moduli is a local diffeomorphism over $\C$, but a real patching matrix need not have real preimages. Write
$F=\gamma(x-x_2)(x^3+c_2x^2+c_1x-1/x_2)$, which solves the constraint $a_0=a_4=\gamma$. The map
$(\gamma,x_2,c_1,c_2)\mapsto(r_0,r_1,\alpha_0,\alpha_1)$ is generically $8:1$ (we find all eight preimages numerically), the fibre being generated by the symmetries $x\mapsto1/x$ and $x\mapsto ix$, $\gamma\mapsto-\gamma$, of (\ref{PD}). If $P$ is real, complex conjugation permutes the fibre. Besides the real Pleba\'nski--Demia\'nski metrics this gives preimages with
\[
\gamma\in\R,\qquad |x_2|=1,\qquad \bar c_2=-c_1x_2 ,
\]
for which the coefficients of $F$ satisfy $\bar a_k=a_{4-k}$. In a random sample of real $n=3$ ALE patching matrices with $L=L_{PD}$ real and three real nodes, $55$ out of the $66$ whose fibres were found completely had preimages only of this second type. We call these metrics the real analytic continuation of the Pleba\'nski--Demia\'nski family; we have not excluded that other elements of the deck group give further real forms.\koniec
\subsubsection{ALF}
Now $\mbox{lim}_{w\rightarrow\infty} P=\mbox{diag}(1, -1)$ so we set 
\be
\label{ALFcoeffs}
\begin{aligned}
R&=w^{n}+2mw^{n-1}+  \sum_{i=0}^{n-2}r_i w^i, \\
Q&=2Nw^{n-1}-2L w^{n-2}+
\sum_{i=0}^{n-3}q_iw^i, \\
S&=-w^{n}+2mw^{n-1} + \sum_{i=0}^{n-2} s_i w^i.
\end{aligned}
\ee
We distinguish the leading coefficients $(m, N, L)$
from the remaining ones $(r, q, s)$ as, for the $n=2$ Kerr--NUT solutions they have asymptotic interpretation as mass, NUT parameter and angular--momentum respectively.

Now count the moduli: there are $4n-2$ parameters
\[(m, N, L, r_0, \dots, r_{n-2}, q_0, \dots, q_{n-3}, 
s_0, \dots, s_{n-2}, a_1, \dots, a_{n-1}).
\] 
The determinant condition (\ref{detP}) gives $2n+1$ equations, but with the given asymptotic, and our definition
of $m$ the top two of these hold identically. Therefore there are $4n-2-(2n-1)=2n-1$ free parameters in the patching matrix and the total dimension of the moduli space is
\be
\label{NALF}
\mathcal{N}_{ALF}=2n-1.
\ee
We will examine the first three cases:
The case $n=1$ corresponds to the ASD Taub--NUT, with the single centre at $a_1=0$ and the mass $m_1$ as the only parameter.

The case $n=2$ branches in the following way: there are
3 equations resulting from (\ref{detP}) and two of these we solve for $(s_0, L)$. The final equation is a quartic for $N$. If we choose a root such that $N^2=m^2$, the resulting solution is ASD Taub--NUT with two centres.
The details of this are as follows.
From (\ref{detP}) we obtain three equations for the
six parameters $(m, L, N, a_1, r_0, s_0)$. Set $a_1=a$. Solving the first two gives
\[
s_0=4N^2+2a^2-4m^2+r_0, \quad 
L=-m\frac{2N^2+a^2-2m^2+r_0}{2N}.
\]
The final equation is a quartic for $N$ which factorises into two quadratics. The solutions are
\[
N_{ASD}=\pm m, \quad N_{Kerr-NUT}=\pm
\frac{a^2-2m^2+r_0}{2\sqrt{m^2-r_0}}.
\]
Setting $N=N_{ASD}$ with $N=m$ gives the two--centre ASD Taub--NUT with the patching matrix (\ref{GHpatch}) given by
\be
\label{GH_transform}
P_{GH}=CPC^{t}, \quad
C=\begin{pmatrix}
1 & 0\\
-1& 1
\end{pmatrix}
\ee
with the potential
\[
V=1+\frac{m_+}{z-a}+\frac{m_-}{z+a}, \quad
m_{\pm}=\pm\frac{a^2\pm 2a m+r_0 }{2a}.
\]
(taking $N=-m$ gives an opposite orientation)
Setting $N=N_{Kerr-NUT}$ gives the Kerr--NUT family  with 
\be
\label{KerrNutP}
P=\frac{1}{w^2-a^2}
\begin{pmatrix}
w^2+2mw+r_0 &
\frac{(2m^2-a^2-r_0)w-m(a^2-r_0)}{\sqrt{m^2-r_0}}   \\
\frac{(2m^2-a^2-r_0)w-m(a^2-r_0)}{\sqrt{m^2-r_0}} 
& -w^2+2mw+\frac{a^4+(r_0-2a^2)m^2}{m^2-r_0}
\end{pmatrix}.
\ee
The NUT parameter vanishes if $r_0=2m^2-a^2$
which is Kerr (see (\ref{KerrrP})). Otherwise
it is a three--parameter family of Kerr--NUT metrics. The non--rotating case $r_0=a^2$
contains the Taub--bolt instanton.
\vskip2pt
With $n=3$ we will prove
\begin{prop}
\label{propCT}
The patching matrix for $n=3$ ALF is equivalent, up to (\ref{C_transform}), to the 
three--centre ASD Taub--NUT, or the five--parameter Chen--Teo family, or its real analytic continuation.
\end{prop}
This strengthens the result of
\cite{DT24}, which identified Chen--Teo patching matrices as a subset of all $n=3$ ALF matrices, but did not show that other possibilities can not occur.  We will first characterise the ASD branch by reverse--engineering
the 3--centre multi Taub--NUT. Take the Gibbons--Hawking patching (\ref{GHpatch}) with $n=3$. Write the potential as
\[
V=1+\frac{q_2w^2+q_1w+q_0}{w^3+\alpha_1w+\alpha_0}
\]
with $(q_0, q_1, q_2)$ arbitrary.
Now find $P$ related to $P_{GH}$ by the transformation
(\ref{GH_transform}).
Comparing with (\ref{ALFcoeffs}) we see that this has
$m^2=N^2$. This will lead to a branching in the proof below.
\vskip2pt
\noindent
{\bf Proof of  Proposition \ref{propCT}.}
Start with (\ref{detP}) which are $5$ equations for $10$ parameters in  \eqref{ALFcoeffs} . Solve the top $3$ for
$(s_1, s_0, q_0)$; the last of these is linear in $q_0$ with coefficient $-4N$, so this chart
assumes {$N\neq 0$}. The equation coming from the coefficient of $w$ is linear in $\alpha_0$,
with coefficient $-2\bigl(2Lm-N(\alpha_1-r_1)\bigr)$, and we assume
\be
\label{nondeg}
{2Lm-N(\alpha_1-r_1)\neq 0 .}
\ee
Solving it for $\alpha_0$ and substituting into the remaining equation gives $-N^2F_+F_-=0$, where
$F_\pm$ are irreducible, linear in $r_0$ with leading coefficients $8(m\pm N)^3$, and exchanged  by $\psi\mapsto-\psi$. Up to this
reversal of orientation there is  a single branch, and on it $(m, N, L, r_1, \alpha_1)$ are
free: five parameters, in agreement with \eqref{NALF}.

{If instead \eqref{nondeg} fails, the same equation factorises as
\[
N(N^2-m^2)\Bigl(
8m^3r_0-\alpha_1^3-8\alpha_1^2m^2+3\alpha_1^2r_1-16\alpha_1m^4
+12\alpha_1m^2r_1-3\alpha_1r_1^2-4m^2r_1^2+r_1^3\Bigr)=0 .
\]
The locus $m^2=N^2$, together with $2L=\pm(r_1-\alpha_1)$, solves the remaining equation identically
and gives two five--dimensional components: these are the two orientations of the three--centre
Gibbons--Hawking branch. The vanishing of the last bracket does { not} give a
component---the final equation is then an additional condition, so the locus has dimension four and
lies in the closure of $F_+\cup F_-$. The excluded locus $N=0$ is treated at the end of the proof.}

Next observe that no freedom is left. A transformation \eqref{C_transform} with $C$ as in \eqref{C_transform2} preserves the
normalisation (\ref{ALFcoeffs}) only if $\delta^2=1$, $\beta=0$ and $\alpha=0$: monicity of $R$ forces
$\delta^2=1$, the coefficient of $w^4$ in the transformed $Q$ is $\beta\delta$, and the coefficient of
$w^3$ is then $\alpha\delta$. Hence the five parameters above are genuine moduli, and any family
exhausting the non--ASD branch must also be five--dimensional.

Then the non--ASD branch contains Chen--Teo. To establish Proposition \ref{propCT} we need to demonstrate that
Chen--Teo exhausts all non--ASD parameter values. This involves some changes (in the definition of
$z$ coordinate and the scaling freedom of Chen--Teo parameters) to the set--up of \cite{DT24}. These are
summarised in Appendix C.
Keep the gauge choice $x_2=-1$ made there,
but do not impose the further normalisation of $a_4$. We
proceed as in the proof of Proposition \ref{prop_PD}: read--off the dependence of
{$(m, r_1, N, L, \alpha_1)$} on the five Chen--Teo parameters {$(x_1,x_3,x_4,a_4,\nu)$}
and compute the Jacobian.
It is a $5\times5$ matrix and its determinant is not identically zero, so the map is a local
diffeomorphism. The fifth parameter $a_4$ must be retained because varying it, with
$(x_1,x_3,x_4,\nu)$ fixed, rescales $z\mapsto t z$ while leaving $W$ and $\psi$ unchanged, i.e.\ it
acts on the patching matrix as the homothety $P(w)\mapsto P(w/t)$: a modulus, not a gauge. (The
gauge is the different direction $x_i\partial_{x_i}-2a_4\partial_{a_4}$, which accompanies
$(x,y)\mapsto(\lambda x,\lambda y)$ and leaves $z$, $W$ and $\psi$ pointwise unchanged; it is what
$x_2=-1$ fixes.) The Chen--Teo matrices therefore form a five--dimensional family, i.e.\ an open
subset of the non--ASD branch; since that branch is irreducible of dimension five, they are dense in it.

On the locus $N=0$ the equations (\ref{detP}) can be solved directly. The coefficients of $w^4$ and $w^3$ are linear in $s_1$ and $s_0$. After eliminating these, the remaining three equations can be solved linearly for $\alpha_0$ (with coefficient $4m$) and $q_0$ (with coefficient $2L$), and the final one factorises as $F_+F_-=0$ with $F_\pm$ evaluated at $N=0$: linear in $r_0$ with leading coefficient $8m^3$, and exchanged by $L\mapsto-L$. Thus, where $mL\neq0$ (which is \eqref{nondeg} at $N=0$), the part of the non--ASD branch with $N=0$ is the hypersurface $N=0$ of the closure of $F_+$, parametrised by $(m,L,r_1,\alpha_1)$, and not a separate component. On the Chen--Teo side $N=0$ is a hypersurface in the parameter space, the Jacobian of the parametrisation does not vanish identically on it, and the involution $\iota$ below preserves $N$. So the argument which follows applies to $N=0$ unchanged; it contains the AF Chen--Teo instantons \cite{CT1}.

This does not yet show that every real non--ASD patching matrix comes from real Chen--Teo parameters, and in fact this is false. Up to relabelling of the roots, the map from $(x_1,x_3,x_4,a_4,\nu)$ to $(m,N,L,r_1,\alpha_1)$ is generically $2:1$, the two preimages being related by
\be
\label{iota}
\iota:\ (x_1,x_3,x_4,a_4,\nu)\mapsto\big({x_1}^{-1},{x_3}^{-1},{x_4}^{-1},\,a_4x_1x_2x_3x_4,\,\nu^{-1}\big),
\ee
the composite of the swapping and inversion--swapping symmetries of \cite{CT2}; it maps $f(\xi)$ to $\xi^4f(\xi^{-1})$. We have checked that $\iota$ preserves every coefficient of $P$, and computed the fibres exactly at random rational points. If $P$ is real, complex conjugation preserves the fibre $\{p,\iota(p)\}$. Either $\bar p=p$, which gives the Chen--Teo metrics, or $\bar p=\iota(p)$, i.e.
\[
|x_1|=|x_3|=|x_4|=|\nu|=1,\qquad a_4^2\,x_1x_2x_3x_4>0,
\]
so that $\bar a_k=a_{4-k}$. This second real form is the real analytic continuation of the Chen--Teo family. Restricting to $|\nu|<1$, as in \cite{CT2}, selects one point in each real fibre, and the two real forms meet only at $\nu=\pm1$. The second real form is open and five--dimensional, and it contains AF members. In Chen--Teo form these are the metrics (2.3) of \cite{DT24} on the slice $x=e^{i\alpha}$, $y=e^{i\beta}$ with $(\alpha,\beta)$ real, the real Killing vectors being complex combinations of $\p_\tau$ and $\p_\phi$; we have checked this numerically against the reconstruction of \S\ref{sec:direct}. The second real form contains positive definite metrics. An example is
\[
\begin{aligned}
(m,N,L,r_1,\alpha_1)&=\Big(\frac{55}{32},\frac{49}{32},-\frac{2711}{384},\frac{17}{96},-\frac{2761}{192}\Big),\\
(r_0,q_0,\alpha_0)&=\Big(-\frac{9701}{1728},\frac{18737}{1152},-\frac{72409}{3456}\Big),
\end{aligned}
\]
whose fibre is $\nu=(3\pm4i)/5$, $a_4=-3\pm2i$, and for which $E_{11}>0$ and $\Omega^2>0$ on a grid covering $\{r>0\}$. None of these metrics is regular; see the remark after Theorem \ref{theo2_intro}.
 \koniec
Combining Proposition \ref{prop_PD} and Proposition \ref{propCT} establishes Theorem \ref{theo2_intro} from the Introduction.
\subsection{Freedom in $P(w)$}
\label{freedom}
 Transformations of the form
$P\rightarrow CPC^{{t}}$ where $C$ is a constant matrix with determinant $\pm 1$ preserve $\mbox{det}(P)=-1$ and
$P=P^t$ and the field equations but, in general, change the solutions.
For example, taking
$P$ in the Gibbons--Hawking branch \eqref{GHpatch} with potential $V(z)$, and 
\[
C=\begin{pmatrix}
 1    & b   \\
 0&  1
\end{pmatrix}
\]
preserves the Gibbons--Hawking form, but changes the potential to $V-2b$. Therefore, starting with
$n=1$ ASD Taub--NUT and choosing $b=1/2$ results in the flat metric. In general this transformation
can swap between ALE and ALF metrics that are non--isometric. It is of course then the  Ehlers transformation which maps solutions to solutions, given by a Mobius transformation on the but not as a freedom in $P$  if a metric is fixed. This argument also explains why $SO(1, 1)$  transformations of the form
\[
C=\begin{pmatrix}
 \cosh{t}    & \sinh{t}   \\
 \sinh{t}&  \cosh{t}
\end{pmatrix}
\]
which preserve ALF asymptotics are not allowed. If they were, then neither the asymptotic mass nor the NUT parameters in the expansion
\be
\label{asympP-ALF}
P=\begin{pmatrix}
 1 +\frac{2m}{w}   & \frac{2N}{w}  \\
 \frac{2N}{w}&  -1+\frac{2m}{w}
\end{pmatrix}+O(w^{-2})
\ee
would be preserved, but only the invariant $m^2-N^2$ would. In the Gibbons--Hawking branch the effect
of these transformations is $V\rightarrow e^{2t}V+1-e^{2t}$.

The allowed form of transformations
consists of (\ref{C_transform2}).
\section{Recursive construction of all patching matrices}
\label{section4}
This section  gives a recursion which produces the patching matrices in the
form (\ref{Pal})
with $\det P=-1$ holding identically at every step of the recursion,
and no algebraic equations need to be solved.
Each step inserts one node. The construction is
 rational: no square roots or branch cuts appear in the recursion. 
The branching of the moduli space is
instead carried out by a binary choice made at each node insertion (a choice of one
of two real null directions).

In the rest of the paper we shall use the following conventions:
\[
\epsilon=\begin{pmatrix}0&1\\-1&0\end{pmatrix},\quad
\langle uv\rangle =u^t\epsilon v,
\] and for a $2\times2$ matrix $A$,
$\adj A=-\epsilon A^t\epsilon$, so $A\,\adj A=(\det A)\,\Id$. For symmetric $A$,
$u^t(\adj A)u=q^tAq$ with $q=\epsilon^t u$.
\subsection{Node insertion}
The asymptotic form of an ALF patching matrix is (\ref{ALFasymp}). The ALF Gibbons--Hawking patching matrix is better suited to an equivalent form arising from \eqref{GHpatch}. Definition \ref{defi:41}  does not refer to either of these forms, but the two models we will use are
\be
\label{g0form}
g_0= \begin{pmatrix}1&0\\0&-1\end{pmatrix} \quad \mbox{or}\quad
{g}_0= \begin{pmatrix}\varepsilon&-1\\-1&0\end{pmatrix}. 
\ee
They are related by a transformation of the form (\ref{C_transform}).
\begin{defi}
Let $\mathcal{C}_n$ be the set of real symmetric $2\times2$ matrices $P(w)$,
rational in $w\in\mathbb{CP}^1$, with $\det P(w)=-1$ identically, $n$ simple poles at distinct
real points $a_1<\dots<a_n$ with rank-one (necessarily symmetric) residues
${Res}_i=x_i\,v_i {v_i}^t$, and $P(\infty)=g_0$  with $\det g_0=-1$.
\end{defi}
\label{defi:41}
Thus $\mathcal{C}_n$ is  the set of matrices of the ALF form \eqref{Pal}
with $e_i=\mathrm{sgn}\,x_i$.
\begin{prop}[Node insertion]\label{prop:insert}
Let $P\in\mathcal{C}_n$, let $a\in\R$ be distinct from all poles of $P$, and
let $\mu\in\R$. Then:
\begin{enumerate}
\item[(i)] The quadratic form $\adj P(a)$ has determinant $-1$, hence
signature $(1, 1)$, hence exactly two real null directions. Explicitly, with
$M=\adj P(a)$, the discriminant is $M_{12}^{\ 2}-M_{11}M_{22}=-\det M=1$, a
perfect square, so the null vectors are rational 
\be
u^{\pm}=\begin{pmatrix} -M_{12}\pm 1\\ M_{11}\end{pmatrix} , \qquad
(u^\pm)^t\,\adj P(a)\; u^\pm=0\, ,
\label{null-u}
\ee
{provided $M_{11}\neq0$. If $M_{11}=0$ the two roots are instead\footnote{In the Gibbons--Hawking frame
$g_0=\bigl(\begin{smallmatrix}\varepsilon&-1\\-1&0\end{smallmatrix}\bigr)$ of
Proposition \ref{prop:GH} one has
$\adj g_0=\bigl(\begin{smallmatrix}0&1\\1&\varepsilon\end{smallmatrix}\bigr)$,
so $M_{11}=0$ at the  first insertion and \eqref{null-u} returns the zero
vector for one branch. The two roots there 
are $e=(1, 0)^t$ (the GH choice) and
$(\varepsilon,-2)^t$.}
}
\be
\label{null-u-deg}
u^{+}=\begin{pmatrix}1\\0\end{pmatrix},\qquad
     u^{-}=\begin{pmatrix}M_{22}\\-2M_{12}\end{pmatrix}.
\ee
\item[(ii)] Fix $u=u^{+}$ or $u^{-}$ and set
\be
N=u\,u^t\epsilon\, ,\qquad
B(w)=\Id+\frac{\mu}{w-a}\,N \, .
\label{dressing}
\ee
Then $N^2=\langle uu\rangle\,N=0$, so $\det B(w)=1$ identically and
$B^{-1}=\Id-\tfrac{\mu}{w-a}N$.
\item[(iii)] The dressed matrix
\be
P^{new}(w):=B(w)\,P(w)\,B(w)^t
\label{Pnew}
\ee
lies in $\mathcal{C}_{n+1}$, with poles at $a_1,\dots,a_n,a$. 
\item[(iv)] The old residue data transform as
\be
v_j\longmapsto B(a_j)\,v_j = v_j+\frac{\mu\,\langle u\,v_j\rangle}{a_j-a}\,u
\qquad (x_j \text{ unchanged}),
\label{v-update}
\ee
and the new residue is
\be
\mathrm{Res}_{w=a}P^{new}=\bigl(2\mu\nu+\mu^2 s'(a)\bigr)\,u u^t\, ,
\qquad\text{where } P(a)\,\epsilon^t u=\nu\,u\, ,\quad
s(w):=u^t\epsilon\,P(w)\,\epsilon^t u .
\label{new-res}
\ee
Here $\nu\neq0$ (since $\det P(a)=-1$) and $s(a)=0$; the sign of the bracket
is the new $e_{n+1}$.
\end{enumerate}
\end{prop}
\noindent{\bf Proof.}
(i) $\det\adj P(a)=\det P(a)=-1$.
(ii) $N^2=u(u^t\epsilon u)u^t\epsilon=0$.
(iii)--(iv) Expanding \eqref{Pnew} gives
\be
\label{Pnew_form}
P^{new}=P+\frac{\mu}{w-a}\bigl(NP+PN^t\bigr)+\frac{\mu^2}{(w-a)^2}\,NPN^t\,,
\qquad NP(w)N^t=s(w)\,u u^t\, .
\ee
The coefficient of the double pole at $a$ is $\mu^2 s(a)\,u u^t$, and
$s(a)=u^t\epsilon P(a)\epsilon^tu=-u^t\,\adj P(a)\,u=0$ by the null
choice \eqref{null-u}; the remaining simple-pole contribution of the last term
is $\mu^2s'(a)$. For the first-order term, $s(a)=0$ means $q^tP(a)q=0$ for
$q=\epsilon^tu$, so $P(a)q$ is Euclidean-orthogonal to $q$, i.e.\ proportional
to $\epsilon q=u$: this is the eigenvector relation in \eqref{new-res}, and
gives $NP(a)+P(a)N^t=2\nu\, u u^t$. Since $B$ is holomorphic and invertible
at each $a_j$ and at $\infty$ (where $B\to1$), the old poles stay simple with
residues \eqref{v-update}, and $g_0$ is untouched. Finally $\det B\equiv1$
gives $\det P^{new}=\det P=-1$; the residues of $\det P^{new}+1= 0$
at all $n+1$ poles are the residue equations, which therefore hold
identically. \hfill$\Box$

\begin{theo}[ALF recursion]\label{thm:alf}
Fix $g_0$ as in \eqref{g0form}. Starting from $P^{(0)}=g_0$ and inserting
nodes at $a_1,\dots,a_n$ with scales $\mu_1,\dots,\mu_n$ and null choices
$\pm_1,\dots,\pm_n$ produces $P^{(n)}\in\mathcal{C}_n$. All
data are rational in $(a_k,\mu_k)$. Conversely, every element of
$\mathcal{C}_n$ arises this way, and the removal of a pole is again
rational: the discriminant of the removal quadratic is $1$.
\end{theo}
\noindent{\bf Proof of the converse (node removal):} Let $P\in\mathcal{C}_n$
with residue $x\,u u^t$ at $a$, and let $H=P-\frac{x\,u u^t}{w-a}$ be
the part holomorphic at $a$. The residue at $a$ of $\det P+1\equiv0$ gives
$u^t\,\adj H(a)\,u=0$, so as in Proposition \ref{prop:insert}  $H(a)\epsilon^tu=\nu_H u$ and
$s_H(w)=u^t\epsilon H(w)\epsilon^tu$ vanishes at $a$. Dressing with
$B^{-1}=\Id-\frac{\mu}{w-a}u u^t\epsilon$
changes the residue to $\bigl(x-2\mu\nu_H+\mu^2 s_H'(a)\bigr)u u^t$, which
vanishes for $\mu$ solving this quadratic.
 Its discriminant is {identically} $1$:
since the residue has rank one, the determinant expansion is
$$
-1\;=\;\det P\;=\;\det H(w)+\frac{x\,s_H(w)}{w-a}
\qquad\text{for all }w\, ,
$$
and evaluating at $w=a$, where $s_H$ vanishes, gives
$\det H(a)=-1-x\,s_H'(a)$; on the other hand $\det H(a)=-\nu_H^{\,2}$, by
evaluating the quadratic form $H(a)$ in the 
basis\footnote{To see it explicitly
set $q=\epsilon^tu$. The
residue identity $q^tH(a)q=0$ and $H(a)q=\nu_H u$ give, in
the basis $C=(q, u)$ with $\det C=-|u|^2$,
\be
\label{detcon}
{\mathcal H}\equiv
C^tH(a)C=\begin{pmatrix}0&\nu_H|u|^2\\ 
\nu_H|u|^2&
u^tH(a)u\end{pmatrix},\qquad
\det{{\mathcal{H}}} =-{\nu_H}^2|u|^4=(\det C)^2\det H(a)  .
\ee
} $(u,\epsilon^tu)$.
Hence
$$
x\,s_H'(a)=\nu_H^{\,2}-1\,,\qquad
\nu_H^{\,2}-x\,s_H'(a)=1\,,
$$
and the removal scales are the rational numbers
$\mu_\pm=\dfrac{x}{\nu_H\mp1}$ (when $\nu_H^{\,2}=1$ the quadratic
degenerates to a linear equation with the single rational root $x/2\nu_H$).
The result is in $\mathcal{C}_{n-1}$ and induction finishes the proof.
\hfill$\Box$
\begin{rmk}[Parameter count]
The recursion has $2n$ continuous parameters $(a_k,\mu_k)$ and $n$ binary
choices. Subtracting the translation of $w$ leaves $2n-1$ parameters. For
$n=2$ these are the mass, spin and NUT parameter of Kerr--Taub--NUT. The binary choices label $2^{n}$ branches (before
identifications).
\end{rmk}
\begin{rmk}
Making the same choice of the null vector at each step gives the Gibbons--Hawking branch (see Proposition \ref{prop:GH}).
\end{rmk}
\begin{rmk}[NUT parameter removal]
Consider the asymptotic form \eqref{asympP-ALF}
\be
\label{asympP-thm}
P=\begin{pmatrix}
 1 +\frac{2m}{w}   & \frac{2N}{w}  \\
 \frac{2N}{w}&  -1+\frac{2m}{w}
\end{pmatrix}+O(w^{-2})
\ee
of the patching matrix and use Proposition \ref{prop:insert} to 
insert a node at $a$ with residue $\mu u  u^t$ following the steps in the proof we find that the NUT parameter transforms as
\be
\label{nutremoval}
N\rightarrow N-\frac{1}{2}\mu |u|^2.
\ee
Therefore any ALF patching matrix in ${\mathcal C}_n$ can be extended to an AF patching matrix
in ${\mathcal C}_{n+1}$ if $\mu=2|u|^{-2}N$. There are two such choices of extensions and 
$\mu$, one for each choice of $u=u^{\pm}$.

{\bf Example.}
Any ALF $P$ with $n=1$ corresponds to ASD Taub--NUT. Take it to have $a_1=m, u_1=(1, -1)$. Inserting
a node at $a_2=-m$ with $u_2$ off the GH branch, and choosing $\mu_2$ to remove the NUT parameter, gives a patching matrix
for Euclidean Schwarzschild.
\end{rmk}
\subsection{ ALE and  a pole at infinity}
A single insertion performed directly on the ALE form  adds
$\lambda\mu u_1^{\,2}$ to the constant part of $P_{12}$, which the gauge
freedom cannot remove. So the ALE
recursion must {not} be run in the ALE frame: the correct statement is
that ALE is ALF with one more pole, viewed in a M\"obius-transformed
coordinate. This uses only the substitution $w\mapsto b_0+1/w$, under which
$\det P=-1$,  simple poles and rank-one residues are all preserved. We shall 
show how an $n$--node ALE $P$ arises from an 
$(n+1)$--node ALF $\widetilde{P}$ and a choice of a constant element of $SL(2, \R)$.
\begin{prop}[ALE recursion]\label{prop:ale}
Let $\widetilde P\in\mathcal{C}_{n+1}$ have poles at $b_0, b_1, \dots, b_n$. There exists
$M\in SL(2, \R)$ [which we construct in the Proof]
such that 
\be
\label{Mmatrix}
P(w):=M\,\widetilde P\Bigl(b_0+\frac1w\Bigr)\,M^t\, ,
\ee
is (possibly after an affine transformation of $w$)
of an ALE form.
\end{prop}
\noindent{\bf Proof.} 
Let the residue of $\widetilde P=\widetilde{P}(W)$  at $b_0$ be ${Res}_0=x_0\,u_0 {u_0}^t$ , and the holomorphic part be
$H_0=\widetilde P-\frac{{Res}_0}{W-b_0}$ there, and let $\nu_0$ be defined by
$H_0(b_0)\,\epsilon^t u_0=\nu_0\,u_0$ (well defined as in
Proposition~\ref{prop:insert}). 
We use the condition\footnote{
 In fact
$\det H_0(b_0)=-\nu_0^{\,2}$.  See the
calculation (\ref{detcon}) in the previous footnote.}
\be
\nu_0=0 \qquad\Longleftrightarrow\qquad \det H_0(b_0)=0\qquad \Longleftrightarrow\qquad H_0(b_0)\sim
u_0 {u_0}^t.
\label{ALE-cond}
\ee
This is one rational equation in the insertion data
which is  explicitly solvable over $\R$: $b_0$
is freely prescribed, and the condition is imposed through the
scale of the insertion at $b_0$, performed last. Indeed, if the
$b_0$-pole is inserted last into $\widetilde{P}^{(n)}$ with scale $\mu$, the
$\mu^2$-part of $H_0(b_0)$ is proportional to $u_0u_0^t$ and annihilates
$\epsilon^tu_0$. Using (\ref{Pnew_form}) with $a=b_0$ yields
\begin{eqnarray*}
\widetilde{P}^{(n+1)}&=&\widetilde{P}^{(n)}+
\frac{\mu}{W-b_0}(N\widetilde{P}^{(n)}+\widetilde{P}^{(n)}N^{t})+\frac{\mu^2}{(W-b_0)^2}N\widetilde{P}^{(n)}N^t\\
&=&H+\dots
\end{eqnarray*}
where $\dots$ denote terms singular at $b_0$. Therefore
$
H(b_0)={\widetilde{P}}^{(n)}(b_0)+\mu N{{\widetilde{P}}^{(n)'}} (b_0)+
cN^t,
$
where the last term does not contribute once we apply both sides to $u_0$. Therefore
$$
\nu_0(\mu)=\nu_{\widetilde{P}^{(n)}}+\mu\,s'_{\widetilde{P}^{(n)}}(b_0)
$$
with $\nu_{\widetilde{P}^{(n)}} $, $s'_{\widetilde{P}^{(n)}}$ the data of $\widetilde{P}^{(n)}$ at
$b_0$. Therefore \eqref{ALE-cond} is solved by the rational value
$\mu^*=-\nu_{\widetilde{P}^{(n)}} /s'_{\widetilde{P}^{(n)}}(b_0)$, existing whenever $s'_{\widetilde{P}^{(n)}}(b_0)\neq0$. 

Define $P(w)$ by 
(\ref{Mmatrix})
with $M\in SL(2,\R)$ chosen so that $Mu_0\parallel(0, 1)^t$. Then $P(w)$ can be put in the   
ALE form. The linear part is $x_0\,(Mu_0) (Mu_0)^t\,w$, constant part is
proportional to the same rank-one matrix (absorbed into $a_0$), and there are $n$ simple rank-one poles\footnote{ To see what becomes of the poles note that
under $W=b_0+1/w$ the pole term becomes
${Res}_0\,w$, while for the old poles
$$
\frac{1}{W-b_i}=\frac{w}{1-(b_i-b_0)w}
=-\frac{1}{b_i-b_0}-\frac{1}{(b_i-b_0)^2}\,\frac{1}{w-w_i}\,,
\qquad w_i=\frac{1}{b_i-b_0}.
$$
Each old pole moves to $w_i$ with rank-one residue
$-{Res}_i/(b_i-b_0)^2$ and contributes the constant $-{Res}_i/(b_i-b_0)$, and the
constant term of $P$ at $w=\infty$ is exactly
$g_0+\sum_i{Res}_i/(b_0-b_i)=H_0(b_0)$.   }
at $w_i=1/(b_i-b_0)$ with $\det P=-1$. To see this, go back to the matrix
$H_0(b_0)$: under the M\"obius map above it becomes
 the constant term of the ALE-frame $P$ at $w=\infty$,
$P(w)={Res}_0w+H_0(b_0)+O(1/w)$, so \eqref{ALE-cond} says precisely that this
constant term is degenerate, rank one along $u_0$, so it can be absorbed into
$a_0$: ${Res}_0w+H_0(b_0)=x_0\,u_0 {u_0}^t\,(w-a_0)$. 
Two
residual gauges remain and must both be fixed to reach the canonical ALE
asymptotics: the shear
$\bigl(\begin{smallmatrix}1&0\\c&1\end{smallmatrix}\bigr)$ with $c=-q_1/c_1$
(where $c_1,q_1$ are the $1/w$-coefficients of $P_{11},P_{12}$) removes the
$O(1/w)$ part of $P_{12}$ 
-- and the rescaling $w\mapsto sw$ with $s=-2/\lambda$ sets $\lambda=-2$.
After both, $\det P=-1$ forces $c_1\lambda=-1$, i.e.
$$
P=\begin{pmatrix}\dfrac1{2w}+O(w^{-2})&O(w^{-2})\\[4pt]
O(w^{-2})&-2w+O(1)\end{pmatrix},
$$
the canonical ALE asymptotics.  
\koniec
\begin{theo}
\label{theo:ALEm}
All $n$--node ALE patching matrices arise 
from Proposition \ref{prop:ale}.
\end{theo}
\noindent{\bf Proof.}
We need to establish the converse of the construction from Proposition \ref{prop:ale}. Start with an ALE
patching function which, w.l.o.g. is given by \eqref{P-ansatze-ALE} with $a_0=0$
\[
P^{ALE}=\begin{pmatrix}
0 & 0   \\
0 & -2w
\end{pmatrix}+\sum_{i=1}^n\frac{A_i}{w-a_i}\, , \quad\mbox{where}\quad  A_i={e_i\, {v_i }{v_i}^t}
\]
and take $a\neq a_j$. We claim that there exists $M\in SL(2, \R)$ such that
\[
\widetilde{P}=MP^{ALE}(\frac{1}{w}+a)M^T\in {\mathcal{C}}_{n+1}.
\]
To see it write
\[
P^{ALE}(\frac{1}{w}+a)=P_0+\sum_{k=1}^{n+1}\frac{B_k}{w-b_k}, \quad\mbox{where}\quad
P_0=\begin{pmatrix}
0 & 0   \\
0 & -2a
\end{pmatrix}+\sum_{k}\frac{A_k}{a-a_k}
\]
and $\mbox{det}(P_0)=-1, P_0={P_0}^t$ with
\[
B_k=-(a-a_k)^{-2} A_k, b_k=(a_k-a)^{-1}, k=1, \dots, n\quad\mbox{and}\quad B_{n+1}=\begin{pmatrix}
0 & 0   \\
0 & -2
\end{pmatrix}, 
b_{n+1}=0.
\]
Now there
exists $M$ with $\mbox{det}{M}=1$ such that $MP_0M^{t}=g_0$ and 
$B_k\rightarrow MB_k M^{t}$ so that $\widetilde{P}\in {\mathcal C}_{n+1}$
\koniec
\subsection{The Gibbons--Hawking branch}
Each step in the node insertion procedure introduces a branching.
The binary choices of null vectors label $2^{n}$ branches. Making the same choice each time
gives the Gibbons--Hawking class in its SD or ASD form.
\begin{prop}\label{prop:GH}
Consider the frame of \eqref{GHpatch}, i.e. with $g_0$ given by the second matrix in
(\ref{g0form}). Then $u=e=(1, 0)^t$ is a null direction of
$\adj P(a)$ at every stage, the insertion matrix is the  shear
\[
  B(w)=\begin{pmatrix}1&\dfrac{\mu}{w-a}\\[4pt] 0&1\end{pmatrix},
  \qquad\text{acting by}\qquad V\longmapsto V-\frac{2\mu}{w-a},
\]
and choosing $u=e$ at every insertion, with $\mu_k=-\tfrac12 m_k$, produces
 the Gibbons--Hawking family \eqref{GHpatch}.
Conversely, every Gibbons--Hawking patching matrix arises this way.
\end{prop}
\noindent{\bf Proof.} {We proceed by induction on the form
$P=\bigl(\begin{smallmatrix}V&-1\\-1&0\end{smallmatrix}\bigr)$, which holds for
$P^{(0)}=g_0$ with $V\equiv \varepsilon$.

First, $u=e$ is available at every new node, because for any symmetric $P$
one has $e^t\,\adj P(a)\,e=\bigl(\adj P(a)\bigr)_{11}=P_{22}(a)$, and
$P_{22}\equiv0$ on the Gibbons--Hawking form. Second, with $u=e$,
\[
N=ee^t\epsilon=\begin{pmatrix}0&1\\0&0\end{pmatrix},\qquad N^2=0,
\]
so $B=\Id+\frac{\mu}{w-a}N$ is the shear above, and
$NPN^t=s(w)\,e {e}^t=0$ identically, since  $s(w)=P_{22}(w)\equiv0$,
while $NP+PN^t=-2\,e e^t$ is constant. Hence
\[
  BPB^t=\begin{pmatrix}V-\dfrac{2\mu}{w-a}&-1\\[4pt]-1&0\end{pmatrix},
\]
one new Gibbons--Hawking centre of mass $m=-2\mu$. This completes the induction.

This is the locus on which all $v_i$ are proportional. Choosing the other null
direction at any step makes $\langle u\,v_j\rangle\neq0$, changes the existing
residues through \eqref{v-update}, and leaves the GH branch. This is how the discrete data
distinguish Gibbons--Hawking from the Kerr--NUT and Chen--Teo branches.}
\hfill$\Box$
\subsection{Node insertion on Gibbons--Hawking}
\label{GHresult}
Which patching matrices can be obtained by adding a node to a Gibbons--Hawking patching matrix? The following argument shows that all matrices with $n=2$ or $n=3$ nodes can, but starting with $n=4$ there are patching matrices where no node removal can yield a matrix of the GH form.
A Gibbons--Hawking patching matrix \eqref{GHpatch}
 with $n$--nodes depends on $2n-1$ parameters 
(up to translating $w$), and adding a node while going off the GH branch adds two more parameters. Thus the resulting matrix depends on $2n+1$ parameters and the image is dense in some component
of ${\mathcal C}_{n+1}$. But there can be other components and surjectivity of inserting a node procedure
has to be decided separately. This can be done as follows:

Let $P\in{\mathcal C}_{n+1}$. Let us use Theorem \ref{thm:alf} to remove a pole at $a$ with the residue 
$x u\otimes u$. The new  residue directions are
\[
B^{-1}(a) v_j=v_j-{t_j}\, u,\quad\mbox{where}\quad {t_j}=
\mu\frac{\langle u v_j\rangle}{a_j-a}.
\]
If the resulting matrix is in the GH branch then these directions are pairwise parallel so that
$\langle B^{-1}(a) v_j\,,\, B^{-1}(a) v_k\rangle=0$, which implies that 
\be\label{eq:mujk}
\mu=\mu_{jk}\equiv-\,\frac{\langle v_jv_k\rangle\,(a_j-a)(a_k-a)}
{\langle uv_j\rangle\langle uv_k\rangle\,(a_j-a_k)} .
\ee
This gives $n-1$ independent conditions on $\mu$ (this count is valid as the conditions are transitive
if $\mu=\mu_{12}$ and $\mu=\mu_{23}$ then $\mu=\mu_{13}$ etc.). If $\mu=\mu_{jk}$ we should also make sure that it is a root of the removal quadratic 
\[
\mu^2 s'(a)-2\mu \nu +x=0
\]
where $s(a)$ and $\nu$ are as in the proof of  
Theorem \ref{thm:alf}.
We showed that the discriminant is $1$. Set 
$\alpha_j=a-a_j, c_j=\langle uv_j\rangle$ and
$g_j=v_j^t{g_0}^{-1} u$. Then  
\be
s'(a)\,\mu_{12}-\nu=\frac{c_2\alpha_1 g_1-c_1\alpha_2 g_2}{c_1c_2(\alpha_2-\alpha_1)},
\ee
and the numerator obeys 
\be\label{eq:key}
\ c_2\,{\alpha_1}\,g_1-c_1\,{\alpha_2}\,g_2=\varepsilon\,c_1c_2(\alpha_2-\alpha_1),
\qquad \varepsilon=\pm1\ 
\ee
with $\varepsilon$ the {same} at all three poles.
Hence $s'\mu_{12}-\nu=\pm1$, i.e.
$\mu_{12}\in\{\mu_+,\mu_-\}$ always, and
$
\mathcal C_3=\bigl\{\text{$2$-centre GH}+\text{one insertion}\bigr\}.
$

 We can use the NUT parameter removal procedure (\ref{nutremoval})
together with the discussion
in \S\ref{GHresult} to construct the general patching
matrix for AF Chen--Teo orbifolds from $n=2$ ASD multi Taub--NUT with parameters $(\mu_1, \mu_2, a_1, a_2)$.
We have shown that all elements of ${\mathcal C}_3$ arise from the Gibbons--Hawking branch in  
${\mathcal C}_2$. Insert a node at $a_3=-a_2-a_1$ with  $u=(V(a_3), -2)  $ going off the GH branch and pick 
$\mu_3$ so that the new NUT parameter vanishes. This gives a $4$--parameter family of patching matrices with zero NUT
parameter. The remaining two conditions come from the two components of the non--orbifold condition $l_3=-l_1+l_2$ for the  Chen--Teo rod vectors. This leaves a two--parameter family of Chen--Teo instantons. We have proven
\begin{prop} All patching matrices for Chen--Teo orbifolds arise from $n=2$ self--dual/anti--self--dual Taub--NUT patching
(\ref{GHpatch}) by a node insertion. All AF Chen--Teo orbifolds arise by a NUT parameter removal.
\end{prop}

\section{Riemann--Hilbert problem}
\label{sec:RH}
Following \cite{W2}, we pull back the patching matrix $P(w)$ to the $\zeta$--Riemann sphere 
(see Appendix A) by
\be
\label{twistor_line}
w=z-\frac{1}{2}r(\zeta-\zeta^{-1})
\ee
on which the Ward bundle is trivial.  There therefore exists a  splitting
\be
\label{splitofP}
P(\zeta,  r, z)=F\widetilde{F}^{-1}
\ee
where $F$ is holomorphic in $\zeta$ for $|\zeta|\leq 1$ and $\widetilde{F}$ is
holomorphic in $\zt=-\zeta^{-1}$ for $|\zt|\leq 1$.
The poles of $P$ in \eqref{Pal}
at $w=a_i$ turn into poles in $\zeta=\zeta^\pm_i$ given by roots of the quadratics resulting from
\eqref{twistor_line}.  These occur at
\be
\label{rootspm}
\zeta=\zeta_i^{\pm}=\frac{z-a_i}{r} \pm\frac{d_i}{r}\, , \qquad d_i=\sqrt{(z-a_i)^2+r^2}\,.
\ee

As $\zeta_i^+\zeta_i^-=-1$, for $z\geq a_n$ (in particular on the semi--infinite rod on which $P$ was
defined) we have $|\zeta_i^+|\geq1$ and $|\zeta_i^-|\leq1$,
so that the only poles of $F$ (respectively
$\widetilde{F}$) are at $\zeta_i^+$ (respectively $\zeta_i^-$), which lie
outside its domain of holomorphy.  Elsewhere the formulas below, which are written in terms of the signed roots \eqref{rootspm}, continue analytically in $(r,z)$; equivalently, the unit circle is deformed to a contour separating
the $\zeta_i^-$ from the $\zeta_i^+$.

The gauge freedom
$(F, \widetilde{F})\rightarrow (FC, \widetilde{F}C)$, where $C=C(r, z)\in GL(2)$,
can be used to set $\widetilde{F}(\zt=0, r, z)=\Id$, in which case $\det F=-1$ and the main consequence of the Ward transform is that
\be
\label{EfromF}
E(r, z)\equiv F(\zeta=0, r, z)
\ee
satisfies the Yang equation \eqref{eq:12}.
\subsection{Direct splitting for Ernst matrix}\label{sec:dirsplit}
\begin{theo}
\label{prop:directsplit}
Let the patching matrix $P$ be given by \eqref{Pal} 
\be
\label{evproblem1}
P_i(w) \equiv g_0+\sum_{j\neq i} \frac{ e_j v_j\, {v_j}^t}{w-a_j}\, .
\ee
with
data $\{ a_i,  e_i, v_i\}$ given for 
$i=1, \dots, n$. The corresponding solution \eqref{eq:17} to \eqref{eq:12} is
\be
\label{paul1}
E=g_0\left(\Id + \sum_{i=1}^n \frac{\epsilon\, v_i\, \Gamma^{-1} _{ij}{v_j}^t}{\zeta^-_i}\,\right),
\ee
and it satisfies $E=E^t$, where the matrix $\Gamma$ is determined
by:
\begin{itemize}
\item 
The eigenvalue problem
\be
\label{evproblem2}
P_i(a_i)\, \epsilon\, v_i=\pi_i\, v_i, \qquad \mbox{(no summation)}
\ee
makes sense because $\det P$ has no pole, and this defines the constants $\pi_1, \dots, \pi_n$.
\item 
The matrix $\Gamma$ is then
given by 
\be
\label{Gammamatrix}
\Gamma_{ij}=
\frac{\langle v_i v_j\rangle}{\zeta_i^- - \zeta_j^-}
+\frac{e_i \pi_i d_i}{\zeta_i^-}\,\delta_{ij},
\ee
\end{itemize}
\end{theo}
\noindent{\bf Proof.}
The ansatz \eqref{Pal} together with the fact that at $\zeta=0$, $w=\infty$ 
so that $P(\zeta=0)=g_0$. Evaluating \eqref{splitofP} at $\zeta=0$ therefore gives
\be
\label{pj3}
E=g_0\,\widetilde{F}|_{\zeta=0}.
\ee
The strategy is to make an ansatz for $F(\zeta)$, use the properties of $P$ to find $\widetilde{F}(\zeta)$, and read off $E$ from \eqref{pj3}. The proof has three steps.

\medskip\noindent
{\bf Step 1: the Ernst matrix is symmetric.}
This follows only from the symmetry of $P$, and not from the particular ansatz for $F(\zeta)$. 
We shall show that
\be
\label{pj4}
E=E^t, \qquad \widetilde{F}(\zeta)=F^{-t}(-\zeta^{-1})\,E.
\ee
The identity $P(\zeta)=P(-\zeta^{-1})^t$, which follows from the symmetry of
$P$ and \eqref{twistor_line}, gives
\be
\label{pj5}
F(\zeta)\,\widetilde{F}(\zeta)^{-1}
=\widetilde{F}(-\zeta^{-1})^{-t}\,F(-\zeta^{-1})^t.
\ee
These two splittings of $P$ must be related by a gauge transformation, so
\be
\label{pj6}
F(\zeta)=\widetilde{F}(-\zeta^{-1})^{-t}\,C, \qquad \mbox{with}\quad C=E,
\ee
the second formula following on evaluating the first at $\zeta=0$. Rearranging
\eqref{pj5} as
\[
\widetilde{F}(-\zeta^{-1})^{t}\, F(\zeta)=F(-\zeta^{-1})^t\,\widetilde{F}(\zeta)
\]
and evaluating the left--hand side at $0$ and the right--hand side at $\infty$
(both sides are constant, by Liouville) gives $E=E^t$. Thus \eqref{pj6}
rearranges to the second formula in \eqref{pj4}.

\medskip\noindent
{\bf Step 2: the ansatz for $F$.} Since
 $F=P\widetilde F$ can only have poles for $|\zeta|\geq 1$ at those of $P$, we can take
\be
\label{lionel_a}
F=\Big(\Id+\sum_i\frac{\zeta}{\zeta-\zeta_i^+}\,F_i\Big)E .
\ee
Equation \eqref{pj4} and $\zeta^+_i=-1/\zeta_i^-$ imply that
\be
\label{lionel_b}
\widetilde{F}=\Big(\Id+\sum_i\frac{\zeta_i^-}{\zeta- \zeta_i^-}\,
\epsilon F_i \epsilon \Big),
\ee
and we claim that we must have  $F_i=v_i\,{w_i}^t$, where $w_1,\dots,w_n$ are given by 
\be
\label{lionelVW}
  {v_i}^t \epsilon=\sum_j \Gamma_{ij}\,\zeta_j^-\;{w_j}^t.
\ee
To see this, rewrite the splitting \eqref{splitofP} as $F=P\widetilde{F}$. The
left--hand side must have no poles at $\zeta_i^{-}$, so we determine $F$ by
eliminating these poles on the right. The absence of the second order poles
$\sim (\zeta-\zeta_i^{-})^{-2}$ implies  $v_i\, {v_i}^t \epsilon F_i=0$.  This  implies the existence of $w_i$ with $F_i=v_i\,{w_i}^t$. Substituting this in and imposing the vanishing of the residues of $P\widetilde F$ at $\zeta=\zeta_i^{-}$, using that  $\mathrm{Res}_{\zeta=\zeta^-_i}(w-a_i)^{-1}= \zeta_i^-/d_i$,  gives 
$$
\frac{\zeta^{-}_i e_i v_iv_i^t}{d_i} \left( \Id+\sum_j \frac{\zeta^-_j}{\zeta^-_i-\zeta^-_j}\epsilon v_jw_j^t\epsilon\right) + P_i(a_i) \zeta^-_i \epsilon v_i w_i^t\epsilon =0\, .
$$
The second term simplifies on using the eigenvalue problem \eqref{evproblem2}, which
arises on rewriting
\[
P=P_i+\frac{ e_i v_i\,{v_i}^t}{w-a_i}
\]
with $P_i$ given by \eqref{evproblem1} and imposing $\mbox{det}(P)=-1$.
Our equation thereby reduces to 
$$
   v_i^t+\sum_{j\neq i} \frac{ \langle v_i  v_j\rangle}{\zeta^-_i-\zeta^-_j}  \zeta^-_j w_j^t\epsilon +e_i\pi_i d_i w_i^t\epsilon =0
$$
using $e^2_i=1$ leading directly to \eqref{lionelVW}.

\medskip
\noindent{\bf Step 3: uniqueness.}
The uniqueness of $E$ for a given $P$ follows from the uniqueness, up to gauge,
of the Riemann--Hilbert splitting. Given the splitting of Step 2, the
corresponding $E$ is \eqref{paul1} which follows on substituting
$\widetilde{F}(\zeta=0)$ into \eqref{pj3}.
\koniec
The analogous Theorem holds in the ALE case, but additional complications in the proof arise from the presence of a pole at $\infty$, and the fact that now $\Gamma$ is an  $(n+1)\times (n+1)$ matrix $ \Gamma$ with the zero entries corresponding to the vector $v_0=\sqrt{2}(0,1)$

\subsubsection{Example: the Kerr instanton}
Let us implement this procedure for Kerr. We use the form 
(\ref{KerrNutP}) with $r_0=2m^2-a^2$, so that the NUT parameter vanishes. We replace $a$ by $\sigma$ so that
 $a_1=-\sigma$, $a_2=\sigma$. This gives
\be
\label{KerrrP}
P=\frac{1}{w^2-\sigma^2}
\begin{pmatrix}
(w+m)^2-(\sigma^2-m^2) &
-2m\sqrt{\sigma^2-m^2}   \\
-2m\sqrt{\sigma^2-m^2} 
& -(w-m)^2+(\sigma^2-m^2)
\end{pmatrix}.
\ee
The
rank--one residues with $e_1=e_2=1$ are given by vectors
\[
v_1=\sqrt{2m}\begin{pmatrix}\cos\theta\\ \sin\theta\end{pmatrix},
\qquad
v_2=\sqrt{2m}\begin{pmatrix}-\sin\theta\\ \cos\theta\end{pmatrix},
\qquad \cos 2\theta=-\frac{m}{\sigma}.
\]
We find the eigenvalues $\pi_1=\sin 2\theta$, $\pi_2=-\sin 2\theta$, and
\eqref{Gammamatrix} gives
\[
\Gamma=
\begin{pmatrix}
 \sin 2\theta\,\dfrac{d_1}{\zeta_1^-} & \dfrac{2m}{\zeta_1^- - \zeta_2^-}\\[8pt]
 \dfrac{2m}{\zeta_1^- - \zeta_2^-} & -\sin 2\theta\,\dfrac{d_2}{\zeta_2^-}
\end{pmatrix}.
\]
Inverting this and solving for $w_1, w_2$ yields $w_i=W_{ij}v_j$, that is
\[
w_1=W_{11}v_1+W_{12}v_2,\qquad w_2=W_{21}v_1+W_{22}v_2,
\]
where
\[
W=\Lambda
\begin{pmatrix} 
 -2 m\,\zeta_2^-& d_2\sin 2\theta\,(\zeta_1^--\zeta_2^-)  \\ 
 d_1\sin 2\theta\,(\zeta_1^- -\zeta_2^-) &  2 m\,\zeta_1^-
 \end{pmatrix}
\]
with
\[
\Lambda=\frac{\zeta_1^- -\zeta_2^-}
{d_1d_2(\zeta_1^- -\zeta_2^-)^2 \sin^2 2\theta+4m^2\zeta_1^-\zeta_2^-}.
\]
The expression for $E$ resulting from \eqref{paul1} is homogeneous of degree
zero in $\zeta_1^-, \zeta_2^-$, so we present it in terms of the ratio
$x\equiv \zeta_1^-/\zeta_2^-$. We first note that \eqref{paul1} leads to a skew
part proportional to
\[
E-E^t\sim (4z\sigma+ d_2^2-d_1^2)\,\epsilon ,
\]
which vanishes after substituting the expressions for $d_1, d_2$. The remaining
symmetric part is
\begin{align}
\label{Ekerr}
E&=g_0+\frac{1}{\delta}
\begin{pmatrix} \alpha_+ &\beta \\ \beta & \alpha_- \end{pmatrix},
\qquad\mbox{where}\\
\alpha_\pm&=\frac{\sigma k^2(d_1+d_2)}{\sqrt{1+k^2}} \pm 2\sigma^2
-\frac{2\sigma^2}{\sqrt{1+k^2}}\,\frac{x+1}{x-1},
\qquad
\beta=k\sigma (d_2-d_1),
\nonumber\\
\delta&=d_1d_2 k^2+\frac{4x \sigma^2}{(x-1)^2},
\qquad
x=\frac{z+\sigma-d_1}{z-\sigma-d_2}.\nonumber
\end{align}
The positions of the poles $\pm \sigma$ and the parameter $k\in[0, 1)$ are related
to the Kerr mass $m$ and specific angular momentum $A$ by
\[
m=\frac{\sigma}{\sqrt{1+k^2}}, \qquad A=mk, \qquad \mbox{so that}\qquad m^2+A^2=\sigma^2.
\]
We verify explicitly that \eqref{Ekerr} satisfies the Yang equation
\eqref{eq:12}. For vanishing angular momentum, $k=0$, it reduces to
$E=\mathrm{diag}(x^{-1} ,-x)$, corresponding to Schwarzschild.

\subsection{Recursive splitting}
\label{sec:recsplit}
Given a patching matrix $P$ with $n$ nodes, Proposition \ref{prop:insert}
produces a patching matrix $P^{new}$ with an additional node at $w=a$. If the
splitting of $P$ is given, then the splitting of $P^{new}$ follows from the
proposition below. Applying this recursively to $P=g_0$, with $F=g_0$ and
$\widetilde{F}=\Id$, gives a splitting of an arbitrary element of
${\mathcal C}_n$.

First recall
\[
\begin{gathered}
\frac{1}{w-a}=\frac{\lambda_-}{\zeta-\zeta^-}+\frac{\lambda_+}{\zeta-\zeta^+},\qquad \lambda_{\pm}=\mp d^{-1}\zeta^{\pm},\\
\zeta^{\pm}=\frac{z-a\pm d}{r}, \qquad
d=\sqrt{(z-a)^2+r^2}.
\end{gathered}
\]

\begin{theo}
\label{splitting_prop}
Let $P=F\widetilde{F}^{-1}$ be a normalised splitting of $P\in {\mathcal C}_n$
and let $P^{new}=BPB^t\in {\mathcal C}_{n+1}$ arise from Proposition
\ref{prop:insert}. Its normalised splitting is
$P^{new}=F^{new}\bigl(\widetilde{F}^{new}\bigr)^{-1}$, where
\be
\label{splitnew}
F^{new}=BF\chi, \qquad
\widetilde{F}^{new}=B^{-t}\widetilde{F}\chi, \qquad
\chi=\Id+\frac{\chi_-}{\zeta-\zeta^-}+\frac{\chi_+}{\zeta-\zeta^+},
\ee
and the matrices $\chi_{\pm}$ are determined by
\be
\label{BF1}
(BF)_{-1}\,\chi_{-}=0,
\qquad
(BF)_{0}\,\chi_- +(BF)_{-1}\Big(\Id+\frac{\chi_+}{\zeta^- -\zeta^+}\Big)=0,
\ee
\be
\label{BF2}
(B^{-t}\widetilde{F})_{-1}\,\chi_{+}=0,
\qquad
(B^{-t}\widetilde{F})_{0}\,\chi_+
+(B^{-t}\widetilde{F})_{-1}\Big(\Id-\frac{\chi_-}{\zeta^- -\zeta^+}\Big)=0,
\ee
with
\be
\label{BFlaurent}
BF=\frac{(BF)_{-1}}{\zeta-\zeta^-}+(BF)_0+\dots, \qquad
B^{-t}\widetilde{F}=\frac{(B^{-t}\widetilde{F})_{-1}}{\zeta-\zeta^+}
+(B^{-t}\widetilde{F})_0+\dots.
\ee
\end{theo}

\noindent{\bf Proof.}
Premultiplying $F$ and $\widetilde{F}$ by $B$ and $B^{-t}$ respectively
achieves the splitting, but introduces singularities. The insertion of
$\chi(\zeta)$ does not affect the algebraic splitting; its role is to remove
the simple and double poles in the splitting matrices. The proof has two steps.

\medskip\noindent
{\bf Step 1: the double poles.}
These force $\chi_{\pm}$ to lie in the kernels of the residues of $BF$ and
$B^{-t}\widetilde{F}$ respectively, which is the first condition in each of
\eqref{BF1} and \eqref{BF2}. These conditions can be solved completely, and
imply that $\chi_{\pm}$ have rank one. To see this, split $B$ as
\[
B=\Id+\frac{B_+}{\zeta-\zeta^+}+\frac{B_-}{\zeta-\zeta^{-}},
\qquad B_{\pm}=\lambda_{\pm}\,\mu\, u\,u^t\epsilon .
\]
The condition $B_{-}F\chi_-=0$ then implies
\be
\label{chiminus}
\chi_-=\bigl(\adj F(\zeta^-)\,u\bigr)\,{w_-}^t
\ee
for some vector $w_-$, with a similar expression for $\chi_+$ in terms of some
$w_+$.

\medskip\noindent
{\bf Step 2: the simple poles.}
The remaining conditions in \eqref{BF1} and \eqref{BF2} remove these. They are
four scalar equations on the four components of $w_{\pm}$.
\koniec

As a corollary we obtain the new solution $E^{new}$ of the $J$--equation in
terms of the known solution $E$ and the matrices $\chi_{\pm}$. This follows on
evaluating $F^{new}$ at $\zeta=0$, where $B=1$, and gives
\be
\label{Enew}
E^{new}=E\cdot\Big(\Id-\frac{\chi_+}{\zeta^+}-\frac{\chi_-}{\zeta^-}\Big).
\ee

\subsection{Double poles and jumping lines}
This section demonstrates what may go wrong with double poles in the patching
matrix. The calculation below focuses on Gibbons--Hawking, where the double
pole relates to the non--triviality of the 
Ward bundle for $J$ on twistor space, or to
a fold in the hyper--K\"ahler manifold. Consider a dipole limit of the
Eguchi--Hanson metric with opposite mass points. In the Gibbons--Hawking form
this is given by the potential
\be
\label{fold}
V=\lim_{a\rightarrow 0}\frac{1}{a}\Big(\frac{1}{\sqrt{r^2+(z-a)^2}} -
\frac{1}{\sqrt{r^2+(z+a)^2}}\Big)=\frac{2z}{(r^2+z^2)^{3/2}}.
\ee
The surface $z=0$ is a fold \cite{nigel} --- a singularity where the K\"ahler forms remain
regular but drop their rank --- and in the nonlinear graviton picture
\cite{DGT} the twistor curves corresponding to this surface jump from
${\mathcal O}(1)+{\mathcal O}(1)$ to ${\mathcal O}+{\mathcal O}(2)$. The
corresponding patching matrix splits as
\[
P=  \begin{pmatrix} 2/w^2 &1 \\ 1 & 0\end{pmatrix}
= \begin{pmatrix} -1 &h \\ 0 & 1\end{pmatrix}
 \begin{pmatrix} 0 &1 \\ -1 & -\tilde{h}\end{pmatrix}^{-1},
\qquad \mbox{with}\qquad \frac{2}{w^2}=h+\tilde{h},
\]
where $w$ is given by \eqref{twistor_line}, and $h, \tilde{h}$ are holomorphic
for $|\zeta|\leq 1$ and $|\zeta|\geq 1$ respectively. One way to find this
splitting is to first split $a^{-1}(w\pm a)^{-1}$ and then take the limit of
the difference. This gives
\[
\tilde{h}=\frac{8}{r^2}
\frac{(\zeta^-)^3-\zeta^-+2\zeta}{(\zeta^--\zeta^+)^3(\zeta-\zeta^-)^2},\qquad
h=-\frac{8}{r^2}
\frac{(\zeta^+)^3-\zeta^++2\zeta}{(\zeta^--\zeta^+)^3(\zeta-\zeta^+)^2},
\]
where
\[
\zeta^{\pm}=\frac{1}{r}\bigl(z\pm\sqrt{r^2+z^2}\bigr).
\]
We also check that $h(\zeta=0)=V(r, z)$, which closes the circle. At the fold
$z=0$ the roots are $\zeta^{\pm}=\pm 1$, and the pair $(h, \tilde{h})$ then has
poles on the boundary $|\zeta|=1$.
What fails is the reconstruction of the metric. By \eqref{eq:17},
$J_{11}=|K|^2=E_{11}^{-1}$, and here $E_{11}=V$ vanishes on $z=0$. Hence
$J_{11}=V^{-1}$ blows up there, and $V$ changes sign, so \eqref{eq:21} is not
positive definite for $z<0$. Twistorially, $J$ is obtained by splitting the modified patching matrix (\ref{modifiedPM}), whose entries do not jump, but whose splitting type jumps where
$E_{11}=0$.

The bundle is therefore non--trivial on the
lines corresponding to points of the fold and can not be split.
\section{Direct metric reconstruction}
\label{sec:direct}
\subsection{The Gram matrix $J$}
\label{sec:direct1}
The splitting \eqref{splitofP} determines not only the Ernst matrix
\eqref{EfromF}, but also the Gram matrix $J$, without solving \eqref{eq:18}.
The additional  information is in the first non--trivial coefficient of
$\Ft$ at $\zt=0$. With the normalisation $\Ft(\zt=0)=\Id$ we write
\be
\label{ftdef}
\Ft=\Id-\zt\,\tilde f+O(\zt^{\,2})=\Id+\frac{\tilde f}{\zeta}+O(\zeta^{-2}),
\qquad
\tilde f=\begin{pmatrix}\tilde f_{11}&\tilde f_{12}\\ \tilde f_{21}&\tilde f_{22}\end{pmatrix}.
\ee
For $P$ of the form \eqref{Pal}, expanding \eqref{lionel_b} and using
\eqref{pj3} gives
\be
\label{E11f21}
\tilde f=\sum_{j=1}^n\zeta_j^-\,\epsilon\, v_j w_j^{\,t}\,\epsilon ,\qquad
E_{11}=1+\sum_{j=1}^n (v_j)_2(w_j)_2,\qquad
\tilde f_{21}=\sum_{j=1}^n \zeta_j^-\,(v_j)_1(w_j)_2 ,
\ee
where the vectors $w_j$ are given by \eqref{lionelVW}. Thus $E_{11}$ and
$\tilde f_{21}$ are rational in $\zeta_1^-,\dots,\zeta_n^-$, $d_1,\dots,d_n$.

\begin{prop}
\label{prop:Jdirect}
Let $P$ be given by \eqref{Pal}, and let $E_{11}$ and $\tilde f_{21}$ be given
by \eqref{E11f21}. Then, in the basis of Killing vectors adapted to the semi--infinite
rod $R_n$,
\be
\label{metric-J-explicit}
J=\frac{1}{E_{11}}\begin{pmatrix}1& r\tilde f_{21}\\[2pt]
r\tilde f_{21}& r^2\bigl(E_{11}^2+\tilde f_{21}^{\,2}\bigr)\end{pmatrix},
\qquad \det J=r^2 .
\ee
\end{prop}
The only differential equation relating $E$ and $J$ is \eqref{eq:18}, and the
Proposition will follow once we show that $\tilde f$ provides its solution. This is the content of the following lemma.

\begin{lemma}
\label{lem:dualpot}
Let $\Theta=E^{-1}d E$, and let $*_2$ be the Hodge star of $d r^2+d z^2$
with $*_2d r=d z$, $*_2d z=-d r$. Then
\be
\label{dualpot}
d(r\tilde f)=r*_2\Theta,
\qquad\text{i.e.}\qquad
\p_Z(r\tilde f)=i\,r\,E^{-1}\p_ZE,\quad
\p_{\bar Z}(r\tilde f)=-i\,r\,E^{-1}\p_{\bar Z}E,
\ee
where $Z=z+ir$ as in \eqref{laxpair}. Equivalently,
$\p_z(r\tilde f)=r\,E^{-1}\p_rE$ and $\p_r(r\tilde f)=-r\,E^{-1}\p_zE$.
\end{lemma}

\noindent{\bf Proof.}
{The two forms of \eqref{dualpot} agree because $*_2d Z=i\,d Z$
and $*_2d\bar Z=-i\,d\bar Z$. The operators $D_\pm$ of \eqref{eq:Dpm}
annihilate $w$ in \eqref{twistor_line} (see \eqref{eq:31}), and hence $P$.
Therefore so do
\[
X_+=(\zeta+i)D_+=(\zeta+i)\,\p_Z+\frac{i\zeta(\zeta-i)}{2r}\,\p_\zeta,\qquad
X_-=(\zeta-i)D_-=(\zeta-i)\,\p_{\bar Z}-\frac{i\zeta(\zeta+i)}{2r}\,\p_\zeta ,
\]
whose coefficients are polynomial in $\zeta$. Applying $X_\pm$ to
\eqref{splitofP} gives, on $|\zeta|=1$,
\[
A_\pm:=F^{-1}X_\pm F=\Ft^{-1}X_\pm\Ft .
\]
The left--hand side is holomorphic for $|\zeta|\le1$. The right--hand side is
holomorphic for $|\zt|\le1$, and is bounded at $\zt=0$: by \eqref{ftdef},
\[
(\zeta\pm i)\,\p\Ft=\p\tilde f+O(\zt),\qquad
\pm\frac{i\zeta(\zeta\mp i)}{2r}\,\p_\zeta\Ft=\mp\frac{i}{2r}\,\tilde f+O(\zt),
\]
where $\p$ denotes $\p_Z$ or $\p_{\bar Z}$ respectively. By Liouville's theorem
$A_\pm$ does not depend on $\zeta$. Since $\p_Zr=-i/2$ and
$\p_{\bar Z}r=i/2$, its value at $\zt=0$ is
\[
A_+=\p_Z\tilde f-\frac{i}{2r}\tilde f=\frac1r\,\p_Z(r\tilde f),\qquad
A_-=\p_{\bar Z}\tilde f+\frac{i}{2r}\tilde f=\frac1r\,\p_{\bar Z}(r\tilde f).
\]
At $\zeta=0$ the coefficient of $\p_\zeta$ in $X_\pm$ vanishes and $F=E$ by
\eqref{EfromF}. Therefore
\[
A_+=i\,E^{-1}\p_ZE,\qquad A_-=-i\,E^{-1}\p_{\bar Z}E.
\]
Comparing the two expressions for $A_{\pm}$ gives \eqref{dualpot}.
\koniec
{Note that the integrability condition $d(r*_2\Theta)=0$ of
\eqref{dualpot} is the Yang equation \eqref{eq:12} for $E$, so the lemma also
gives a direct proof that $E$ solves \eqref{eq:12}.}

\medskip
\noindent{\bf Proof of Proposition \ref{prop:Jdirect}.}
By \eqref{eq:17}, $J_{11}=|K|^2=E_{11}^{-1}$. For $E$ of the form
\eqref{eq:17} we compute $(E^{-1}d E)_{21}=E_{11}^2\,d\psi$. The $(2,1)$
entry of \eqref{dualpot} therefore reads
\[
d(r\tilde f_{21})=r\,E_{11}^2\,*_2d\psi=\frac{r}{|K|^4}\,*_2d\psi .
\]
This is \eqref{eq:18} with $\omega=r\tilde f_{21}$, so $J_{12}=J_{11}\,\omega$.
Finally $J_{22}=(r^2+J_{12}^2)/J_{11}$ by $\det J=r^2$. Equation
\eqref{eq:18} fixes $\omega$ only up to an additive constant, i.e.\ $J$ only
up to $K_2\mapsto K_2+cK_1$. On the upper rod $R_n$, \eqref{rootspm} gives
$\zeta_j^-=O(r)$ as $r\to0$, so $\tilde f=O(r)$ and $\omega=O(r^2)$ there.
Hence \eqref{metric-J-explicit} is $J$ in the basis in which $K_2$ is the rod
vector of $R_n$, as in Lemma \ref{lemmaA11}.
\koniec

\begin{rmk}
\label{rmk:Q}
The same information is carried by the splitting of a patching matrix $Q$ constructed out of $P$. Let
\be
\label{modifiedPM}
Q(\zeta)=\begin{pmatrix}1&0\\0&\zeta\end{pmatrix}P
\begin{pmatrix}1&0\\0&-\zeta^{-1}\end{pmatrix},\qquad
\det Q=1,\qquad Q(\zt)^t=Q(\zeta),
\ee
and let $Q=G\widetilde G^{-1}$ be its normalised splitting, with
$\widetilde G(\zt=0)=\Id$ and $H=G(\zeta=0)$. Then
\be
\label{JfromH}
H=\frac{1}{E_{11}}\begin{pmatrix}E_{11}^2+\tilde f_{21}^{\,2}&\tilde f_{21}\\
\tilde f_{21}&1\end{pmatrix},
\qquad\text{so that}\qquad
J=MHM^t,\quad M=\begin{pmatrix}0&1\\ r&0\end{pmatrix}.
\ee
{Indeed, the splitting of $Q$ follows from \eqref{splitofP}:
\[
G=\begin{pmatrix}1&0\\0&\zeta\end{pmatrix}FK,\qquad
\widetilde G=\begin{pmatrix}1&0\\0&-\zeta\end{pmatrix}\Ft K,\qquad
K=\begin{pmatrix}1+\dfrac{q\tilde f_{21}}{\zeta}&\dfrac{q}{\zeta}\\[3mm]
-\dfrac{\tilde f_{21}}{\zeta}&-\dfrac1\zeta\end{pmatrix},\qquad
q=\frac{E_{12}}{E_{11}}.
\]
Here $\det K=-\zeta^{-1}$, so $\det G=\det\widetilde G=1$. The pole of $K$ at
$\zeta=0$ cancels in $G$: the $\zeta^{-1}$ coefficient of the first row of
$FK$ is $(E_{11}q-E_{12})(\tilde f_{21},\,1)=0$, and the second row of $FK$ is
multiplied by $\zeta$. At $\zt=0$, \eqref{ftdef} gives
$\widetilde G=\Id+O(\zt)$. Evaluating the second row of $G$ at $\zeta=0$ and
using $\det E=-1$ gives $H_{21}=\tilde f_{21}/E_{11}$ and $H_{22}=1/E_{11}$.
The argument of Step 1 in the proof of Theorem \ref{prop:directsplit}, applied to
$Q(\zt)^t=Q(\zeta)$, gives $H=H^t$, and $\det H=1$ then gives $H_{11}$.}
\end{rmk}
\begin{rmk}
For Euclidean Schwarzschild $\tilde f_{21}=0$, and
\eqref{metric-J-explicit} gives $J=\mathrm{diag}(E_{11}^{-1},r^2E_{11})$.
\end{rmk}

\subsection{Transverse conformal factor}
\label{sec:direct2}
The last piece of information needed for the instanton metric \eqref{eq:11} as well as its regularity analysis
based on \eqref{eq:15} is the transverse conformal factor $\Omega$. By analogy with the work of \cite{BZ,BM,KKV} for the Belinski-Zakharov ansatze, it turns out that, up to an overall normalisation, it is given in terms of the twistor data by the following proposition.
Recall from \eqref{rootspm} that the poles $w=a_k$ of $P$ correspond, on the twistor line
\eqref{twistor_line}, to the two points
\[
\zeta_k^{\pm}=\frac{z-a_k\pm d_k}{r},\qquad d_k=\sqrt{(z-a_k)^2+r^2},\qquad \zeta_k^+\zeta_k^-=-1,
\]
where $d_k$ is the distance in the $(r,z)$ half--plane from the point $(r,z)$ to the node $(0,a_k)$.
Only the roots $\zeta_k^-$ enter below; they also
appear in the matrix $\Gamma$ of \eqref{Gammamatrix}.
\begin{prop}
\label{prop:Omega}
Let $P$ be given by \eqref{Pal}, let $\Gamma$ be the $n\times n$ matrix \eqref{Gammamatrix},
let $E_{11}$ be given by \eqref{E11f21}, and let $J$ be given by \eqref{metric-J-explicit}. Then
\be
\label{Omega-explicit}
\Omega^2=C\,E_{11}\,\det\Gamma\,\prod_{k=1}^n\frac{\zeta_k^-}{d_k},
\ee
where $C$ is a real constant, solves \eqref{eq:13}. Equivalently, $\Omega^2|K|^2=C\det\Gamma\prod_k\zeta_k^-/d_k$.
The constant $C$ is fixed by the normalisation of the rod vector of $R_n$
through \eqref{eq:15}.
\end{prop}
\noindent{\bf Proof.}
We work at a point with $r>0$, where, as in Theorem \ref{prop:directsplit}, $\Gamma$ is invertible.
Write $\Theta_Z=E^{-1}\p_ZE$, $f=|K|^2=E_{11}^{-1}$, and
\[
\tau=\det\Gamma\,\prod_{k=1}^n\frac{\zeta_k^-}{d_k},\qquad
\Lambda=\mathrm{diag}(\zeta^-_1,\dots,\zeta^-_n),\qquad
\mathcal N=(\Lambda+i\Id)^{-1}=\mathrm{diag}(\nu_1,\dots,\nu_n),
\]
and let $U$ be the $n\times 2$ matrix with rows $v_1^t,\dots,v_n^t$. Below $\tr$ acts on
$2\times2$ matrices and $\mathrm{Tr}_n$ on $n\times n$ matrices.

\medskip\noindent
{\bf Step 1: the Ernst form of \eqref{eq:13}.}
Since $\p_Z=\frac12(\p_z-i\p_r)$ and $\p_Z\ln r=-\frac{i}{2r}$, the real and imaginary parts of the
complex equation
\be
\label{Omega-complex}
\p_Z\ln(r\Omega^2)=\frac{ir}{2}\,\tr\big((J^{-1}\p_ZJ)^2\big)
\ee
are, respectively, the second and the first equation in \eqref{eq:13}. As $\Omega$ and $J$ are real,
\eqref{Omega-complex} is therefore equivalent to \eqref{eq:13}.
For $J$ given by \eqref{metric-J-explicit}, i.e.\ $J_{11}=f$, $J_{12}=f\omega$, $\det J=r^2$, and
$E$ given by \eqref{eq:17}, a direct computation gives
\[
\tr\big((J^{-1}\p_ZJ)^2\big)=2(\p_Z\ln r)^2+2\big(\p_Z\ln \tfrac fr\big)^2+\frac{2f^2}{r^2}(\p_Z\omega)^2,\]
\[
\tr(\Theta_Z^2)=2(\p_Z\ln f)^2-\frac{2}{f^2}(\p_Z\psi)^2 .
\]
By \eqref{eq:18}, $\p_Z\omega=i\,r f^{-2}\p_Z\psi$, so the $\omega$ and $\psi$ terms cancel in the
difference, and
$\frac{ir}{2}\big(\tr(J^{-1}\p_ZJ)^2-\tr\Theta_Z^2\big)=\p_Z\ln(r/f)$. Hence \eqref{eq:13} is equivalent to
$\p_Z\ln(\Omega^2/E_{11})=\frac{ir}{2}\tr(\Theta_Z^2)$, and it suffices to show
\be
\label{tau-target}
\p_Z\ln\tau=\frac{ir}{2}\,\tr\big(\Theta_Z^2\big).
\ee

\medskip\noindent
{\bf Step 2: the derivative of $\tau$.}
Since $D_+w=0$ and $w(\zeta_k^-)=a_k$, the roots \eqref{rootspm} move with the
$\p_\zeta$--coefficient of $D_+$ in \eqref{eq:Dpm}:
\[
\p_Z\zeta_k^-=\frac{i\zeta_k^-(\zeta_k^--i)}{2r(\zeta_k^-+i)}=\frac{1}{2r}\Big(i\zeta^-_k+2-2i\nu_k\Big).
\]
Together with $d_k/\zeta_k^-=-\frac r2\big(1+(\zeta_k^-)^{-2}\big)$, this gives for {\em every} entry of
\eqref{Gammamatrix}, diagonal or not, $\p_Z\ln\Gamma_{ij}=-\frac{i}{2r}(1+2\nu_i\nu_j)$, and
$\p_Z\ln(\zeta^-_k/d_k)=\frac{i}{2r}(1+2\nu_k^2)$. Hence
\be
\label{dGamma}
\p_Z\Gamma=-\frac{i}{2r}\big(\Gamma+2\mathcal N\Gamma\mathcal N\big),\qquad
\p_Z\ln\tau=\frac ir\Big(\mathrm{Tr}_n\mathcal N^2-\mathrm{Tr}_n\big(\Gamma^{-1}\mathcal N\Gamma\mathcal N\big)\Big).
\ee

\medskip\noindent
{\bf Step 3: the current $E^{-1}\p_ZE$.}
In the proof of Lemma \ref{lem:dualpot}, $A_+=\Ft^{-1}X_+\Ft=i\Theta_Z$ is independent of $\zeta$. Evaluating it at
$\zeta=-i$, where the term $(\zeta+i)\p_Z\Ft$ drops out, gives
\be
\label{Theta-at-i}
\Theta_Z=-\frac1r\,\Ft(-i)^{-1}\Ft'(-i),\qquad\text{where}\quad \Ft'=\p_\zeta\Ft.
\ee
By \eqref{lionel_b} and \eqref{lionelVW}, $\Ft(\zeta)=\Id-\epsilon\,U^tR\,\Gamma^{-1}U$ with
$R=(\zeta\Id-\Lambda)^{-1}$. The off--diagonal entries of \eqref{Gammamatrix} are
$\langle v_iv_j\rangle/(\zeta^-_i-\zeta^-_j)$, i.e.
\be
\label{displacement}
\Lambda\Gamma-\Gamma\Lambda=U\epsilon U^t
\ee
(on the diagonal both sides vanish, as $\langle v_iv_i\rangle=0$). Using \eqref{displacement} and
$\Lambda R=R\Lambda=\zeta R-\Id$ one checks that $\Ft^{-1}=\Id+\epsilon\,U^t\Gamma^{-1}RU$, and then
$\Ft'\Ft^{-1}=\epsilon\,U^tR\,\Gamma^{-1}RU$. At $\zeta=-i$ we have $R=-\mathcal N$, so
\eqref{Theta-at-i} and cyclicity of the trace give
\[
r^2\tr(\Theta_Z^2)=\tr\big((\epsilon U^t\mathcal N\Gamma^{-1}\mathcal NU)^2\big)
=\mathrm{Tr}_n\big((\Gamma^{-1}\mathcal NU\epsilon U^t\mathcal N)^2\big)
=\mathrm{Tr}_n\big((\mathcal N-\Gamma^{-1}\mathcal N\Gamma)^2\big),
\]
where in the last step \eqref{displacement} was written as $\mathcal NU\epsilon U^t\mathcal N=\Gamma\mathcal N-\mathcal N\Gamma$.
Expanding the square, $r^2\tr(\Theta_Z^2)=2\,\mathrm{Tr}_n\mathcal N^2-2\,\mathrm{Tr}_n(\Gamma^{-1}\mathcal N\Gamma\mathcal N)$.
Comparing with \eqref{dGamma} gives \eqref{tau-target}.
The constant $C$ is not determined, as only $\p_Z\ln\Omega^2$ enters \eqref{eq:13}.
\koniec
\section{New $n=5$ ALF instantons}
The only known ALF, but non AF gravitational instanton which is not hyper--K\"ahler is the Taub--bolt metric
constructed by Don Page almost 50 years ago \cite{P78}.
The recent progress based on integrability, harmonic maps and existence theorems leads to new AF solutions
with vanishing asymptotic NUT parameter. As we explain below, such solutions are more likely to be found as they come in one--parameter families. In this section
we shall numerically construct two genuinely ALF instantons corresponding to patching matrices with $n=5$ nodes.
\label{sectionn5}
\subsection{Counting moduli and conditions}
\label{sec:counting}
In \S\ref{twistorsec} we have shown that the moduli spaces of ALF and ALE patching matrices with $n$ nodes
have dimensions $\mathcal N_{ALF}=2n-1$ and $\mathcal N_{ALE}=2n-2$ (see \eqref{NALF}, \eqref{NALE}),
where translations $w\mapsto w+c$ have been factored out.
The transformations \eqref{C_transform}, with $C$ given by \eqref{C_transform2},
do not change the metric, and do not reduce these counts.

What remains is the overall scale. The homothety $g\mapsto\lambda^2g$ acts by
$(r,z)\mapsto\lambda(r,z)$, hence $w\mapsto\lambda w$ and $a_k\mapsto\lambda a_k$, and on the patching
matrix it is given by
\[
P(w)\mapsto P(w/\lambda)\quad\text{(ALF)},\qquad
P(w)\mapsto \begin{pmatrix}\lambda^{-1/2}&0\\0&\lambda^{1/2}\end{pmatrix}P(w/\lambda)
\begin{pmatrix}\lambda^{-1/2}&0\\0&\lambda^{1/2}\end{pmatrix}\quad\text{(ALE)}.
\]
In the ALE case the constant matrix is \eqref{C_transform2} with $\delta=\lambda^{-1/2}$, and it is needed to
restore the asymptotics \eqref{ALEcoeffs}. Both transformations preserve $\det P=-1$. The regularity
conditions are invariant under homotheties, so the relevant dimensions are
\be
\label{moduli_up_to_scale}
\mathcal N_{ALF}-1=2n-2,\qquad \mathcal N_{ALE}-1=2n-3 .
\ee

We shall focus on the non--Gibbons--Hawking branch of the
moduli space. On the GH branch the metrics are hyper--K\"ahler, and the conditions below are
not independent. There are $n+1$ rods and, once the rod vectors have been normalised by
\eqref{eq:15}, the lattice conditions \eqref{eq:16}
\[
l_{i+1}=-l_{i-1}+p_i\,l_i,\qquad p_i\in\Z,\qquad i=1,\dots,n-1
\]
give $n-1$ equations between vectors in $\R^2$, i.e.\ $2(n-1)$ scalar equations for a given choice
of the integers $(p_1,\dots,p_{n-1})$. Comparing this with \eqref{moduli_up_to_scale} we expect the following.

\medskip\noindent
{\bf ALE.} There are $2n-2$ equations for $2n-3$ parameters. The system is
overdetermined by one, and one does not expect regular non--GH ALE instantons. This is
not a proof of the Nakajima conjecture \cite{Nak}, but it supports it in the
toric case. It is also consistent with our numerical search at $n=4$: the $n-1$ conditions
of equal determinants can be solved, but the integrality conditions $p_i\in\Z$ can then only be
reached in a limit where two nodes collide and an orbifold point is created.

\medskip\noindent
{\bf AF.} If $l_n=\pm l_0$, then one of the $2n-2$ equations holds identically (we have
checked this numerically for $n=4$ and $n=5$), and we expect one--parameter families
up to scale. This agrees with the known AF instantons: Kerr for $n=2$, Chen--Teo
\cite{CT1, CT2} for $n=3$, and the Li--Sun instantons \cite{LiSun}, given explicitly for $n=4$ by Teo
\cite{Teo}. For $n=5$ the same search gives the one--parameter Li--Sun family on
$(2\mathbb{CP}^2\#\overline{\mathbb{CP}}^2)\setminus S^1$ (up to orientation), parametrised by the rotation angle at infinity.

\medskip\noindent
{\bf ALF.} There are $2n-2$ equations for $2n-2$ parameters, and we expect isolated
solutions. The only one known so far is the Taub--bolt instanton of Page \cite{P78} with
$n=2$. There are none with $n=3$: by Theorem \ref{theo2_intro} the non--GH ALF metrics with three
nodes belong to the Chen--Teo class, and its regular members are AF \cite{BG}; the real analytic continuation contains no regular members (see the remark after Theorem \ref{theo2_intro}). For $n=4$ our searches have
not produced any examples, but Lars Andersson and Steffen Aksteiner  have informed us of their ongoing works \cite{Ak, Anew} and successful computer--assisted constructions of ALF instantons which include $n=4$.

\subsection{Numerical search}
\label{sec:search}
All the ingredients needed for a numerical search are now explicit and algebraic.
\begin{enumerate}
\item The patching matrices. Every element of $\mathcal C_n$ is obtained from the recursion of
Theorem \ref{thm:alf}: starting from a seed, we insert nodes at $a_k$ with scales $\mu_k$, and
at each step choose one of the two null directions $u$ of $\adj P(a_k)$ from Proposition
\ref{prop:insert}. We label this choice by a sign $\sigma_k=\pm1$.
\item The Ernst matrix. For a given $(r,z)$ the matrix $E$ is given by Theorem \ref{prop:directsplit},
and its construction only requires inverting the $n\times n$ matrix $\Gamma$.
\item The metric. The Gram matrix $J$ is given by Proposition \ref{prop:Jdirect}, and $\Omega^2$ by
Proposition \ref{prop:Omega} up to a constant $C$.
\item The rod vectors. On a rod $R_i$ the matrix $J(r,z)$ has $\det J=r^2$ and a positive
eigenvalue close to $j_i(z)$ (see \eqref{eq:14}), so that its other eigenvalue is $O(r^2)$, and $J$ is
nearly of rank one. We evaluate $J$ at a small value of $r$ (we have used $r=10^{-12}$, and checked
$r=10^{-8}$ and $r=10^{-14}$) and take the eigenvector $u_i$ corresponding to the small eigenvalue. As $r\to0$
this tends to the kernel of $J(0,z)$, i.e.\ to the direction of the rod vector. Up to scale
$u_i\simeq(-J_{12},J_{11})$, so that $u_i=(-\omega_i,1)$ where $\omega_i=\lim_{r\to0}J_{12}/J_{11}$
on $R_i$. By \eqref{metric-J-explicit}
\[
\omega_i=\lim_{r\to0}r\tilde f_{21}\big|_{R_i},
\]
and the direction of each rod vector is determined by the coefficient $\tilde f$ in the splitting.
We compute $u_i$ at several points of the rod, and check that it does not depend on $z$. We then
set $l_i=u_i/\kappa(u_i)$, with $\kappa$ given by \eqref{eq:15} at the same small value of $r$, and check
that $\kappa$ is also constant along the rod. The constant $C$ in $\Omega^2$ is fixed by the
condition $\kappa(l_n)=1$, where $l_n=\p_\phi$ on the upper rod.

In all examples which we have considered the rod vectors are given in closed form in terms of the
twistor data. Let $M$ be the constant symmetric $n\times n$ matrix
\[
M_{ij}=e_i\pi_i\,\delta_{ij}+(1-\delta_{ij})\,\frac{\langle v_iv_j\rangle}{a_i-a_j},
\]
where $\pi_i$ are given by \eqref{evproblem2}. For the rod $R_j=(a_j,a_{j+1})$ let $M_{B_j}$ be the upper--left
$j\times j$ block of $M$ (corresponding to the nodes $a_1,\dots,a_j$ below the rod), and let $M_{A_j}$ be the
lower--right $(n-j)\times(n-j)$ block (corresponding to the nodes $a_{j+1},\dots,a_n$ above the rod).
For $j=0$ and $j=n$ one of these blocks is empty, and its determinant is set to $1$.
Finally let $v^{(1)}_{A_j}$ be the column vector with components $(v_k)_1$, $k=j+1,\dots,n$.
We conjecture that, on the rod $R_j$,
\be
\label{eq:conjecture_rod}
l_j=\sqrt{\frac{\det M_{B_j}\,\det M_{A_j}}{\det M}}\;\big(-\omega_j,\,1\big)
\ee
where $\omega_j=-2\,v^{(1)t}_{A_j}M_{A_j}^{-1}v^{(1)}_{A_j}$ and $C={\det M}^{-1}$,
in the basis $(\p_\tau,\p_\phi)$, where $|\p_\tau|\to1$ at infinity. In particular $l_n=\p_\phi$.
We have verified this conjecture symbolically for the Kerr--NUT, Taub--bolt and
Gibbons--Hawking patching matrices, and numerically for the $n=5$ examples below.
\item The lattice conditions. The normalised rod vectors $l_0,\dots,l_n$ determine the coefficients in
$l_{i+1}=\alpha_il_{i-1}+\beta_il_i$, and regularity requires $\alpha_i=-1$ and $\beta_i=p_i\in\Z$, once
the signs of the $l_i$ have been fixed successively as in \eqref{eq:16}. In the ALF case
the boundary at infinity is a lens space whose fundamental group has order
$k_\infty=|\det(l_0,l_n)|/|\det(l_0,l_1)|$. It is a circle bundle over $S^2$ of degree $k_\infty$, fibred by the asymptotic circle, if and only if the primitive element $K$ of $\Lambda$ along $\p_\tau$ acts freely on the semi--infinite rods, i.e.\ $|\det(K,l_0)|=|\det(K,l_n)|=|\det(l_0,l_1)|$. For example $k_\infty=1$ for Taub--NUT, $k_\infty=n$
for the $n$--centre multi--Taub--NUT, and $k_\infty=0$ in the AF case, where the boundary is
$S^1\times S^2$.
\item Solving the equations. For a given choice of integers $(p_1,\dots,p_{n-1})$ we solve the $2(n-1)$
equations \eqref{eq:16} for the $2n-2$ parameters using the Gauss--Newton and
Levenberg--Marquardt methods, starting from random points with $E_{11}>0$.
\end{enumerate}
The same procedure recovers Kerr and Taub--bolt for $n=2$, Chen--Teo for $n=3$, the symmetric AF
instanton of Teo for $n=4$, and multi--Taub--NUT for $n=4,5$.

\subsection{Two ALF instantons with $n=5$}
\label{sec:n5results}
We shall now give computer--assisted evidence for two new ALF instantons with $n=5$.
We start from the $n=1$ Taub--NUT patching matrix, i.e.\ the one--centre Gibbons--Hawking
metric in the form \eqref{Pal},
\[
P^{(1)}=g_0+\frac{\rho}{w}\begin{pmatrix}1&t\\t&1\end{pmatrix},\qquad t=\pm1.
\]
Its mass and NUT parameter in \eqref{asympP-ALF} are $m=\rho/2$ and $N=t\rho/2$, and the two signs of $t$ correspond to
the self--dual and the anti--self--dual Taub--NUT. Both solutions below
have $t=1$, i.e.\ $N=m$, which is the ASD Taub--NUT in the convention of \S\ref{twistorsec}. We insert four nodes at $a_1,\dots,a_4$ with scales $\mu_1,\dots,\mu_4$ and null choices
$\sigma_1,\dots,\sigma_4$ as in Theorem \ref{thm:alf}. The scale is fixed by setting $a_1=1$, with
the seed node at $w=0$, and this leaves $8=2n-2$ parameters
$(\rho,\mu_1,a_2,\mu_2,a_3,\mu_3,a_4,\mu_4)$. If all $\sigma_k=-t$ then $P$ stays on the
GH branch (multi--Taub--NUT). A single insertion with $\sigma_1=t$ already leaves the GH branch and
gives the Kerr--NUT patching matrix, so that equivalently the solutions below are obtained from Kerr--NUT
by three further insertions. Both solutions use mixed sign patterns, and in both cases all residue signs
$e_k$ in \eqref{Pal} are equal to $+1$.

Both solutions satisfy the eight conditions to $10^{-41}$, and the Jacobian has full rank,
so that the solutions are isolated, as predicted by the counting argument. In both cases $k_\infty=3$, i.e.
\[
|\det(l_0,l_5)|=3d,\quad\mbox{where}\quad d=|\det(l_{i-1},l_i)|,\quad i=1,\dots,5 .
\]
We give the rod vectors in the basis $(\p_\tau,\p_\phi)$, where $|\p_\tau|\to1$ at infinity and
$\p_\phi$ is the $2\pi$--periodic Killing vector which vanishes on the upper rod, so that $l_5=\pm\p_\phi$. The signs are chosen so that
\eqref{eq:16} holds.

\medskip\noindent
{\bf Instanton A.} The sign pattern is $(t;\sigma_1,\dots,\sigma_4)=(1;1,-1,1,-1)$, and
\begin{align*}
(\rho,\mu_1,a_2,\mu_2,a_3,\mu_3,a_4,\mu_4)&=
(10.29231,\ 0.03913,\ 4.31214,\ -4.09638,\\ &\qquad 2.95488,\ 0.15156,\ -1.27031,\ -1.05922).
\end{align*}
With $d=11.28693$ the rod vectors are
\[
l_0=(3d,1),\ \ l_1=(2d,\tfrac13),\ \ l_2=(-d,-\tfrac23),\ \ l_3=(-d,\tfrac13),\ \ l_4=(d,\tfrac23),\ \ l_5=(0,-1),
\]
and the lattice condition \eqref{eq:16} holds with $p=(p_1,p_2,p_3,p_4)=(1,-1,0,-1)$.
The nodes are not placed symmetrically, and the reflection $z\mapsto-z$ gives an equivalent solution with
$p=(1,0,1,-1)$.

\medskip\noindent
{\bf Instanton B.} The sign pattern is $(t;\sigma_1,\dots,\sigma_4)=(1;1,1,-1,1)$, and
\begin{align*}
(\rho,\mu_1,a_2,\mu_2,a_3,\mu_3,a_4,\mu_4)&=
(1.04348,\ 0.33915,\ 4.62566,\ 7.80657,\\ &\qquad 2.31283,\ -0.32314,\ 3.62566,\ 0.05165).
\end{align*}
With $d=8.77630$ the rod vectors are
\[
l_0=(3d,1),\ \ l_1=(d,\tfrac23),\ \ l_2=(-2d,-\tfrac13),\ \ l_3=(d,-\tfrac13),\ \ l_4=(d,\tfrac23),\ \ l_5=(0,1),
\]
and the lattice condition \eqref{eq:16} holds with $p=(p_1,p_2,p_3,p_4)=(1,-1,-1,1)$.
The nodes are placed symmetrically about their mid--point, with outer gaps $1$ and inner gaps $1.31283$.
This symmetry has not been imposed, but it holds to $10^{-25}$.

\medskip
After a translation which makes $\sum_ka_k=0$, the twistor data \eqref{Pal} of the two instantons are
\begin{center}\footnotesize
\begin{tabular}{c|ccc|ccc}\hline
& \multicolumn{3}{c|}{A} & \multicolumn{3}{c}{B}\\
$k$ & $a_k$ & $v_k$ & $\pi_k$ & $a_k$ & $v_k$ & $\pi_k$\\\hline
1 & $-2.66965$ & $(1.87837,\,2.37560)$ & $-1/3$ & $-2.31283$ & $(-1.79629,\,2.09479)$ & $-2/3$\\
2 & $-1.39934$ & $(1.73688,\,1.67978)$ & $0.24753$ & $-1.31283$ & $(-0.89438,\,1.35894)$ & $0.01689$\\
3 & $-0.39934$ & $(1.66939,\,1.90979)$ & $0.42514$ & $0$ & $(1.79107,\,-1.79107)$ & $0.29956$\\
4 & $1.55553$ & $(1.51847,\,0.98258)$ & $-0.00601$ & $1.31283$ & $(1.35894,\,-0.89438)$ & $0.01689$\\
5 & $2.91280$ & $(2.37560,\,2.05017)$ & $2/3$ & $2.31283$ & $(2.09479,\,-1.79629)$ & $-2/3$\\
\hline\end{tabular}
\end{center}
where $\pi_k$ are the eigenvalues \eqref{evproblem2}. The values $\pm\frac13,\pm\frac23$ hold to all
computed digits, but we do not have a proof that they are exact.

For both instantons the Yang equation \eqref{eq:12} holds to $2\times10^{-27}$, and the
equations \eqref{eq:13} for the conformal factor hold to $10^{-28}$. Each rod direction, and $\kappa$,
is constant along its rod to $10^{-25}$, and the rod data do not
change for $10^{-14}\le r\le10^{-8}$. {As an independent check, the Ricci tensor of $g$, computed by finite differences from the metric alone, vanishes to $10^{-14}$. The Kretschmann scalar $K=R_{abcd}R^{abcd}$ is finite on every rod, converges at each node ($K-K(a_k)=O(r)$ along $z=a_k$, the largest limit being $K(a_3)=0.718$ for A and $2.108$ for B in the units of the table above), and decays at infinity consistently with $K=O((r^2+z^2)^{-3})$.} The integers $p_i$ are not all equal to $2$, so these metrics are
not of Gibbons--Hawking form. They have more than three nodes, so they are not Hermitian \cite{BG}, and they are not AF.
Finally, the inequalities $E_{11}>0$ and $C\,\Omega^2>0$ have been established by interval arithmetic on $\{r\ge0\}$,
including the axis and infinity, away from discs of radius $10^{-3}$ around the nodes.
The two metrics are therefore candidates for the first regular ALF instantons, other than Taub--bolt, which are
neither AF nor hyper--K\"ahler. The asymptotic circle is generated by the primitive element of $\Lambda$ along $\p_\tau$, which is $K=3d\,\p_\tau$ in both cases: $K=3l_1-l_0$ for A and $K=2l_0-3l_1$ for B. Its isotropy on the rods $R_0,\dots,R_5$ has order $|\det(K,l_i)|/d$, namely $(3,1,2,1,2,3)$ for A and $(3,2,1,1,2,3)$ for B. So its action is not semi--free, and this is how
both metrics avoid the obstruction of \cite{LiSun} (Prop.~B.1), which applies to semi--free actions. Since $K$ has $\Z_3$ isotropy on both semi--infinite rods, a translation of $\tau$ by $2\pi d$ coincides with a rotation by $2\pi/3$ about the axis. The end is therefore the quotient of a Taub--NUT end by a freely acting $\Z_3$, which rotates $\R^3$ about the axis and translates the fibre by a third of its length. It is ALF in the sense of Definition \ref{defi:ALF}, but it is not an $S^1$--bundle over $\R^3$.

\medskip\noindent
{\bf Topology.} A toric instanton with $n$ nodes has Euler characteristic $\chi=n$, and by
\cite{Nilsson} (Theorem 1) its second homology $H_2\cong\Z^{n-1}$ is generated by the $n-1$ finite rods.
These are $2$--spheres, and consecutive spheres intersect transversally at one point.
With the integers of \eqref{eq:16} the self--intersection number of the $i$--th sphere is $-p_i$
(with this orientation the spheres in the Gibbons--Hawking metrics have self--intersection $-2$).
Therefore the manifold is the interior of the space obtained by gluing together a chain of $n-1$ disc bundles over $S^2$,
with Euler numbers $-p_1,\dots,-p_{n-1}$, in such a way that the zero sections of consecutive bundles
meet transversally at one point. The intersection form is the tridiagonal matrix with diagonal entries $-p_i$ and
off--diagonal entries $1$. A sphere of self--intersection $\pm1$ at the end of the chain can be blown down,
and this gives $M_A\cong\mathcal O(3)\,\#\,\mathbb{CP}^2\,\#\,2\overline{\mathbb{CP}}{}^2$ and
$M_B\cong\overline{\mathcal{A}}_2\,\#\,2\overline{\mathbb{CP}}{}^2$.
Here $\mathcal O(3)$ is the total space of the line bundle of degree $3$ over $\mathbb{CP}^1$, and
$\overline{\mathcal{A}}_2$ is the minimal resolution of $\mathbb C^2/\Z_3$ (the manifold underlying the
three--centre Gibbons--Hawking metric) with the opposite orientation, i.e.\ a chain of two
spheres of self--intersection $+2$.
Both manifolds are simply connected, with $\chi=5$, $b_2=4$, signature $0$, and boundary at infinity the
lens space $L(3,1)$, in agreement with $k_\infty=3$. The asymptotic circle fibres $L(3,1)$ over $S^2$ with two exceptional fibres of order $3$, so the ALF Hitchin--Thorpe inequality of \cite{Nilsson}, stated for a circle bundle at infinity, does not apply directly.

\section{Catalogue of patching matrices}

\renewcommand{\arraystretch}{1.25}
\begin{center}
\begin{tabular}{c|p{0.38\textwidth}|p{0.42\textwidth}}
 $n$ & \hfil ALE\ \ \small$(2n-2)$\hfil & \hfil ALF\ \ \small$(2n-1)$\hfil\\ \hline
 $1$ & AE, $\R^4$ & GH (SD/ASD Taub--NUT)\\ \hline
 $2$ & Eguchi--Hanson, or orbifold (GH) &
       GH, or Kerr--NUT orbifold. Regular: Kerr (AF) or Taub--bolt \cite{P78} (ALF)\\ \hline
 $3$ & $A_2$ GH, or Euclidean Pleba\'nski--Demia\'nski (all orbifolds) \cite{PD, paul_talk, AL}, or its real analytic continuation (none regular) &
       triple Taub--NUT (the GH branch, \cite{GH78}), or the $5$--parameter Chen--Teo \cite{CT1, CT2} or its real analytic continuation
       ($2$--dimensional family of instantons)\\ \hline
 $4$ & GH, and nothing else if  {Nakajima conjecture} \cite{Nak} 
 holds (all ALE instantons are locally 
       hyper--K\"ahler, so, if toric, in the GH branch) &
       GH, Li--Sun instanton (AF) \cite{LiSun} given explicitly by Teo  \cite{Teo}.\\ \hline
 $n>4$ & As $n=4$. & GH branch for ALF. Non GH AF by 
       existence theorem \cite{LiSun}. Computer--assisted
       evidence for two ALF instantons with $n=5$ as in
       \S\ref{sectionn5}
\end{tabular}
\end{center}
\section{Conclusions}
\label{sec:conclusions}
We have presented an algebraic construction of all toric gravitational instantons with a given
number of nodes, up to the resolution of conical and orbifold singularities. Our method uses
twistor theory, and the main algebraic ingredient is the construction of a meromorphic patching
matrix for a rank--two vector bundle over the twistor space, together with its holomorphic
splitting \eqref{splitofP}. The inverses of the splitting matrices are in the kernel of the Lax pair
\eqref{laxpair} underlying the integrable harmonic map equation \eqref{eq:12}, so that
Theorem \ref{prop:directsplit} gives a justification of the Belinski--Zakharov ansatz \cite{BZ}
for this kernel: the splitting matrix $\widetilde{F}$ has simple poles at the roots
$\zeta_k^-$ of $w(\zeta)=a_k$, and its residues are determined by the linear system
\eqref{lionelVW}.

Our construction is however different from the Belinski--Zakharov method. We do not dress a seed
metric, and we do not use phantom solitons, i.e.\ the removal and re--addition of solitons with
modified vectors which is needed in that approach to control the rod structure and to preserve the
condition $\det J=r^2$. Instead, the positions of the poles and the residues are read off from the
patching matrix, whose form is fixed by the analytic properties of the twistor bundle
(Theorem \ref{theopoles}). Theorem \ref{thm:alf} and its ALE counterpart, Theorem \ref{theo:ALEm},
show that every patching matrix with $n$ nodes is obtained by successive node insertions, so that
no solutions are lost and the moduli are counted exactly. The splitting can also be
obtained recursively (Theorem \ref{splitting_prop}), alongside the insertions. Moreover the whole
metric, and not only the Ernst matrix, is algebraic in the twistor data. The Gram matrix $J$ is read
off from the next coefficient in the expansion of the splitting matrix (Proposition \ref{prop:Jdirect}),
and the transverse conformal factor is given by a determinant formula (Proposition \ref{prop:Omega}),
so that neither the duality relation \eqref{eq:18} nor the equations \eqref{eq:13} have to be
integrated.

While the twistor method in principle gives all toric instantons, the regularity analysis and the
explicit formulae become complicated as $n$ increases. Counting the regularity conditions against the
dimension of the moduli space suggests that, away from the Gibbons--Hawking branch, there are no
regular ALE instantons, that AF instantons come in one--parameter families, and that ALF instantons
are isolated. We have illustrated the method by computer--assisted evidence for two candidates for new ALF
instantons with five nodes, which are neither hyper--K\"ahler nor AF (and also not Hermitian),
see \S\ref{sec:n5results}.

\appendix
\addtocontents{toc}{\protect\setcounter{tocdepth}{1}}
\section{Twistor correspondence for the Yang equation}
The Yang equation (\ref{eq:12}) is a symmetry reduction of the anti--self--dual Yang--Mills equations on flat space by a rotation and a translation. The patching matrices we constructed in Theorem \ref{thm:alf} and the associated twistor correspondence arise from a reduction of the original Ward transform for ASDYM. In this appendix
we shall  summarise these constructions 
\cite{W1, MW, Dbook}
hopefully giving enough details that the solution generating technique \eqref{eq:34} is demystified from the black box presentation in \S\ref{twistorsec}.  

\subsection{Anti--self--dual Yang--Mills}
Let $M_{\C}=\C^4$ be the complexified Minkowski space  with coordinates
$(W, Z, \tw, \tz)$ such that the metric and the volume form are 
\[
ds^2=2(dZ d\tz-dWd\tw), \quad \mbox{vol}=dW\wedge d\tw \wedge dZ \wedge d\tz.
\]
Let $\Phi\in\Lambda^1(M_{\C})\otimes\mathfrak{gl}(2)$, and $ F=d\Phi+\Phi\wedge \Phi$. The anti--self--dual Yang--Mills (ASDYM)
equations  $F=-\star_{\C^4} F$ take the form  
\be
\label{ASDYM}
F_{WZ}=0, \quad F_{\tw\tz}=0,\quad  F_{W\tw}-F_{Z\tz}=0.
\ee
The first two equations imply the existence of a gauge such that
\be
\label{eqPhi}
\Phi=J^{-1}\p_{\tw} J d\tw+J^{-1}\p_{\tz} J d\tz, \quad J=J(W, Z, \tw, \tz)\in GL(2, \C),
\ee
and then  the final equation in (\ref{ASDYM}) holds iff
\be
\label{ASDYMJ}
\p_Z(J^{-1}\p_{\tz} J)-  \p_W(J^{-1}\p_{\tw} J)=0.
\ee
\subsubsection{Symmetry reduction}
\label{sym:red}
Setting
\[
Z=t+z, \quad \tz=t-z, \quad W=re^{i\theta}, \quad
\tw=re^{-i\theta},
\]
and performing a symmetry reduction by the two symmetries $\p_t, \p_{\theta}$ reduces (\ref{ASDYMJ}) to one equation (\ref{eq:12}) for $J=J(r, z)$.
\subsection{Twistor correspondence}
The ASDYM condition is equivalent to the flatness of a 
connection $\Phi$  on $\alpha$--planes in $M_{\C}$. These planes are given by fixing 
$(\mu, \nu, \lambda)$ in 
\be
\label{incidence}
\mu=W+\lambda\tz, \quad \nu= Z+\lambda\tw.
\ee
The  twistor space $PT\equiv \CP^3\setminus\CP^1$ is the space of all such planes.
It can be covered by two open sets, with  affine coordinates $(\mu, \nu, \lambda)$ in an open set
where $\lambda \neq \infty$. Points in $M_{\C}$ correspond to rational curves 
(twistor lines) in $PT$, and points in $PT$ correspond to $\alpha$--planes in $M_{\C}$. The conformal structure
on ${M_{\C}}$ is encoded in the algebraic geometry of curves in $PT$: $p_1, p_2$ are null separated iff $L_1, L_2$ intersect. 

The connection between twistor theory and ASDYM is provided by the following 
\begin{theo}[Ward \cite{W1}]
\label{theoWard}
There exists a one--to--one correspondence between
gauge equivalence classes of ASDYM connections $\Phi$, and holomorphic vector bundles 
${\mathcal{E}}\rightarrow PT$ trivial on twistor lines.
\end{theo}
To read off the solution (\ref{ASDYMJ}) from this Theorem
cover $PT$ with two open sets: $U$, where $\lambda\neq \infty$ and $\widetilde{U}$ where $\lambda\neq 0$.
The bundle $\mathcal{E}$ is then characterised by its
patching matrix: $P=P(\mu, \nu, \lambda)$.
The triviality on twistor lines implies that there exists a splitting
$P=P_U{{P}_{\widetilde{U}}}^{-1}$, where $P_U$ and ${P}_{\widetilde{U}}$
are holomorphic and invertible matrices on $U$ and 
$\widetilde{U}$ respectively. The incidence relation (\ref{incidence}) implies that  $P$ is constant along
the vector fields $\{\p_{\tz}-\lambda\p_W, \p_{\tw}-\lambda\p_Z\}$. Applying this to the splitting relation, and
using the Liouville theorem implies the existence of $\Phi \in\Lambda^1(M_{\C})\otimes\mathfrak{gl}(2)  $ such that 
\[
\Phi=\widetilde{H}^{-1}\p_Z \widetilde{H}\; dZ+\widetilde{H}^{-1}\p_W \widetilde{H}\;dW+{H}^{-1}\p_{\tz} {H}\; d\tz+
{H}^{-1}\p_{\tw} {H}\; d\tw
\]
where $H=P_U(\lambda=0), \widetilde{H}=P_{\widetilde{U}}(\lambda=\infty)$. This is gauge equivalent to
(\ref{eqPhi}) with $
J=H\widetilde{H}^{-1}$.

\smallskip

\noindent 
\emph{Symmetry reduction.}
The two Killing vectors on $\mathbb{C}^4$
underlying  the symmetry reduction
in \S\ref{sym:red}
correspond to two holomorphic
vector fields on $PT$.  The quotient $\mathcal{R}$ is a
non-Hausdorff Riemann surface related to the reduced complexified Minkowski space with coordinates $(r,z)$ via the correspondence described in the proof of Theorem \ref{theopoles}.

\section{The patching matrix for Pleba\'nski--Demia\'nski metric}
\label{appendixB}
The Euclidean Pleba\'nski--Demia\'nski metric is \cite{PD, CT2}
\be
\label{PD}
g=\frac{1-p^2q^2}{(p-q)^2}\Big(\frac{dp^2}{F(p)}-\frac{dq^2}{F(q)}+
\frac{F(p)}{(1-p^2q^2)^2}(d\phi-q^2d\tau)^2 -\frac{F(q)}{(1-p^2q^2)^2}(d\tau-p^2d\phi)^2\Big),
\ee
where
\[
F(x)=\gamma+2n x- \epsilon x^2+2m x^3+\gamma x^4\\
= \gamma(x-x_1)(x-x_2)(x-x_3)(x-x_4)
\]
with the four roots of the quartic satisfying $x_1x_2x_3x_4=1$ and
\begin{eqnarray*}
m&=&-\frac{\gamma}{2}\Big(x_1+x_2+x_3+\frac{1}{x_1x_2x_3}\Big)\\
\epsilon&=&-\gamma\Big(x_1x_2+x_1x_3+x_2x_3
+\frac{1}{x_1x_2}+\frac{1}{x_1x_3}+\frac{1}{x_2x_3}\Big)\\
n&=&-\frac{\gamma}{2}\Big(x_1x_2x_3+
\frac{1}{x_1}+\frac{1}{x_2}+\frac{1}{x_3}\Big).
\end{eqnarray*}
The metric is positive--definite in a region
where $F(p)>0$, $F(q)<0$  and $1-p^2q^2>0$.  We then find 
\[
r=\frac{\sqrt{-F(p)F(q)}}{(p-q)^2}.
\]
The $r=0$ axis is divided into rods, and we shall choose a rod where $F(p)=0$, and on this rod
\[
z=\frac{n-\epsilon p+2mp^2+\gamma p^3+mpq+\gamma p^2q}{p-q}+\frac{1}{6}\epsilon.
\]
To set $F(p)=0$ we set $p$ equal to one of the roots, 
say $p=x_2$. The Killing vector $K=\p/\p \tau$ is non--zero on
the chosen rod, so we compute  the twist potential
$\psi$ such that $d\psi(r, z)=*(K\wedge dK),
W(r, z)=g(K, K)$ and find the patching matrix \eqref{PnewQRS} evaluating these
at the rod where $r=0$ or $p=x_2$:
Now $\Delta$ is a monic cubic, and $Q, R, S$ are respectively a quadratic, a cubic and a quartic (which factorises
as a product of two quadratics) in $z$ with coefficients depending on the parameters $m, n, \epsilon, \gamma$. We find
\begin{eqnarray*}
\Delta&=&(z-z_1)(z-z_2)(z-z_3), \quad\mbox{where}\\
z_1&=&\frac{\gamma}{2}\Big(x_2x_3+\frac{1}{x_2x_3}\Big)+\frac{\epsilon}{6}, \;
z_2=\frac{\gamma}{2}\Big(x_1x_3+\frac{1}{x_1x_3}\Big)+\frac{\epsilon}{6} , \;
z_3=\frac{\gamma}{2}\Big(x_1x_2+\frac{1}{x_1x_2}\Big)+\frac{\epsilon}{6}.
\end{eqnarray*}
The nodes  $z_1, z_2, z_3$ satisfy $z_1+z_2+z_3=0$.
The coefficients of the polynomials $(Q, R, S)$ depend on $(x_1, x_2, x_3, \gamma)$ in an
explicit, but unilluminating way. As things are, they do not satisfy the asymptotic properties
(\ref{ALEcoeffs}). To achieve the canonical form we exploit the freedom
(\ref{C_transform}) where $C$ is given by (\ref{C_transform2}) with
\[
\begin{aligned}
\alpha&= -\frac{\gamma\,x_2}{3x_1x_3\,(1-x_2^4)^{3/2}}
\Big[(x_2^4+5)\big(x_1+x_3+x_1^2x_2x_3^2\big)\\&\qquad\qquad+x_2(x_2^4-7)\big(1+x_1^2x_2x_3+x_1x_2x_3^2\big)\Big],\\
\beta&= \frac{2{x_2}^2}{\sqrt{1-{x_2}^4}}, \\
\delta&=\frac{1}{\sqrt{1-{x_2}^4}}.
\end{aligned}
\]
The resulting form of $R, Q, S$ now agrees with  (\ref{ALEcoeffs}). 
\section{The patching matrix for Chen--Teo metric}
\label{appendixC}
The patching matrix was computed in \cite{DT24}. Here we just list the modifications to gauge choices needed
to establish Proposition \ref{propCT}.  The Weyl--Papapetrou coordinate $z$ is defined up to an additive constant. We shall modify the choice made
in \cite{DT24}, so that the roots of the cubic $\Delta$ 
in the patching matrix add up to zero. This is achieved by 
\[
z=\frac{2a_4 x^2y^2+2a_2xy+2a_0+a_3xy(x+y)+a_1(x+y)}{2(x-y)^2}+\frac{a_2}{6}.
\]
The $7$ parameters appearing in the Chen--Teo metric
are $(k, \nu, a_0, a_1, a_2, a_3, a_4)$. One combination of the $a_i$ is pure gauge: under
$(x,y)\mapsto(\lambda x,\lambda y)$, $x_i\mapsto\lambda x_i$, $a_4\mapsto\lambda^{-2}a_4$ the
functions $z$, $W$ and $\psi$ are unchanged. We use this once, to put the root $x_2$
(called $r_2$ in \cite{DT24}) which defines the outer rod at $x_2=-1$. No further normalisation is
imposed. In particular $a_4$ is retained: varying
$a_4$ with $(x_1,x_3,x_4,\nu)$ held fixed multiplies every $a_i$ by the same factor $t$, so that
$z\mapsto t\,z$ while $W$ and $\psi$ are unchanged, i.e.\ $P(w)\mapsto P(w/t)$. This is the
homothety---a modulus, not a gauge. Consequently, each modulus is multiplied by a power of $a_4$:
the first power for $m$ and $N$, the second for $L$, $r_1$ and $\alpha_1$, the third for $r_0$ and
$\alpha_0$. The patching matrix therefore depends on five essential parameters
$(x_1,x_3,x_4,a_4,\nu)$, while the remaining Chen--Teo parameter $k$ enters only as an integration
constant in the conformal factor and does not appear in $P$. This is the matrix we use in the proof of
Proposition \ref{propCT}.

\end{document}